\documentclass[a4paper, notitlepage, 12pt]{article}
\usepackage[total={160mm,250mm}, top=25mm, left=25mm, includefoot]{geometry}

\usepackage{ucs}              % Unicode support in LaTeX
\usepackage[utf8x]{inputenc} 	% We use UTF-8 encoded files only 
\usepackage[T1]{fontenc}     % for the right font encoding
\usepackage{lmodern}         % fonts for good pdf files

\usepackage[singlespacing]{setspace}

\usepackage{color}
\usepackage{amsmath,amsthm,amssymb,amscd,amsbsy}
\usepackage{upgreek}

\usepackage{tensor}

\usepackage[nointegrals]{wasysym}
\usepackage{float}
\usepackage{tabularx,paralist}
\usepackage{graphicx}
\usepackage[breaklinks]{hyperref}
\usepackage{authblk}

\usepackage{multirow} 
\usepackage{array}
\usepackage{caption}

\usepackage{rcs} \RCS $Revision: 1.3.1.1 $ \RCS $Date: 2018/07/06 01:25:26 $ \RCS $Author: hgm $ \RCS $RCSfile: 18_RV-algebra-model.tex,v $ \RCS $Id: 18_RV-algebra-model.tex,v 1.3.1.1 2018/07/06 01:25:26 hgm Exp $

\usepackage[british]{babel}

\usepackage{refdef}
\usepackage{ifontdef}

\newcommand{\ignore}[1]{}

\newcommand{\Lp}{\mrm{L}}

\newcommand{\Hp}{\mrm{H}}
\newcommand{\Ck}{\mrm{C}}

\def\frks#1#2{{\raise0.8ex\hbox{{\leavevmode$\scriptstyle #1$}}}{\kern-0.4ex\hbox{/}}{\raise-0.2ex\hbox{\kern-0.4ex\hbox{{\leavevmode$\scriptstyle #2$}}}}}

\newcommand{\vrho}{\varrho}

\newcommand{\vsigma}{\varsigma}

\newcommand{\vphi}{\varphi}
\newcommand{\vpi}{\varpi}
\newcommand{\vkappa}{\varkappa}

\newcommand{\vek}[1]{\mathchoice{\displaystyle\boldsymbol{#1}}
{\textstyle\boldsymbol{#1}}{\scriptstyle\boldsymbol{#1}}
{\scriptscriptstyle\boldsymbol{#1}}}

\newcommand{\mat}[1]{\mathchoice{\displaystyle\mathbf{#1}}
{\textstyle\mathbf{#1}}{\scriptstyle\mathbf{#1}}
{\scriptscriptstyle\mathbf{#1}}}

\newcommand{\opb}[1]{\vek{{\mathsf{#1}}}}

\newcommand{\ops}[1]{\mathchoice{\displaystyle\mathsf{#1}}
{\textstyle\mathsf{#1}}{\scriptstyle\mathsf{#1}}
{\scriptscriptstyle\mathsf{#1}}}

\newcommand{\tnb}[1]{\mathchoice{\displaystyle\mathboldsans{#1}}
{\textstyle\mathboldsans{#1}}{\scriptstyle\mathboldsans{#1}}
{\scriptscriptstyle\mathboldsans{#1}}}

\newcommand{\tns}[1]{\mathchoice{\displaystyle\mathsans{#1}}
{\textstyle\mathsans{#1}}{\scriptstyle\mathsans{#1}}
{\scriptscriptstyle\mathsans{#1}}}

\newcommand{\EXP}[1]{\mathbb{E}\left(#1\right)}

\newcommand{\im}{\mathop{\mathrm{im}}\nolimits}
\newcommand{\dom}{\mathop{\mathrm{dom}}\nolimits}

\newcommand{\tr}{\mathop{\mathrm{tr}}\nolimits}

\newcommand{\spn}{\mathop{\mathrm{span}}\nolimits}

\newcommand{\di}{\mathrm{d}}

\newcommand{\ii}{\mathchoice{\displaystyle\mathrm i}
{\textstyle\mathrm i}{\scriptstyle\mathrm i}
{\scriptscriptstyle\mathrm i}}

\newcommand{\KL}{Karhunen-Lo\`eve}
\newcommand{\ip}[2]{\langle #1, #2 \rangle}
\newcommand{\bkt}[2]{\langle #1 | #2 \rangle}

\newcommand{\ns}[1]{| #1 |}
\newcommand{\nd}[1]{\| #1 \|}

\newcommand{\var}{\text{var}}
\newcommand{\cov}{\text{cov}}

\definecolor{myred}{rgb}{1, 0.2, 0.2}

\newcommand{\autheadcr}{\authorcr}
\newcommand{\citep}[1]{\cite{#1}}

\newcommand{\authorhgm}{Hermann G. Matthies}

\newcommand{\affilwire}{Institute of Scientific Computing\autheadcr
                        Technische Universit\"at Braunschweig\autheadcr
                        38092 Braunschweig, Germany\autheadcr
                        e-mail: \ttt{wire@tu-bs.de}}

\newcommand{\thetitle}{Analysis of Probabilistic and\\ Parametric Reduced Order Models}

\newcommand{\theauthor}{\authorhgm}
\newcommand{\thesubject}{35B30, 37M99, 41A05, 41A45, 41A63, 60G20, 60G60, 65J99, 93A30}
\newcommand{\thekeywords}{stochastic models, parametric models,
                          correlation, factorisation, spectral decomposition}

\newcommand{\textdate}{\today}

\newcommand{\thebib}{./bib}

\begin{document}

% ============================================================================
% connect to default LaTeX values
\title{\thetitle\thanks{Partly supported by the Deutsche
          Forschungsgemeinschaft (DFG) through SPP 1886 and SFB 880.}}

%\title{\thetitle}

\author{\authorhgm}
%\author[o]{\authorro}

%\makeatletter
%\author{\authorhgm 
%%  Acknowledgements
%\thanks{Corresponding author: D-38092 Braunschweig, 
%       Germany, e-mail: \texttt{wire@tu-bs.de}}
%}
%\makeatother
%\author{\theauthor}

\affil{\affilwire}
%\affil[o]{\affilcnam}

%\date{\today}
%\date{}
\date{\textdate}

%%%% END IGNORE
%%-%%

\maketitle

% !TEX encoding = UTF-8 Unicode
% !TEX root = ../18_RV-algebra-model.tex
% RCSID:       $Id: abstract.tex,v 1.2 2018/07/05 22:01:22 hgm Exp $
% Author:      $Author: hgm $
% Contact:     wire@tu-bs.de
% ============================================================================
%% texfile{
%%  AUTHOR    = "$Author: hgm $",
%%  VERSION   = "$Revision: 1.2 $",
%%  DATE      = "$Date: 2018/07/05 22:01:22 $",
%%  FILENAME  = "$RCSfile: abstract.tex,v $"}
%
% =================================

\begin{abstract}
Stochastic models share many characteristics with generic parametric
models.  In some ways they can be regarded as a special case.  But
for stochastic models there is a notion of weak distribution or
generalised random variable, and the same arguments can be used
to analyse parametric models.
Such models in vector spaces are connected to a linear map, and in
infinite dimensional spaces are a true generalisation.
Reproducing kernel Hilbert space and affine- / linear- representations in terms
of tensor products are directly related to this linear operator.  This
linear map leads to a generalised correlation operator, and representations
are connected with factorisations of the correlation operator.
The fitting counterpart in the stochastic domain to make this point
of view as simple as possible are algebras of random variables with
a distinguished linear functional, the state, which is interpreted as expectation.
The connections of factorisations of the generalised correlation
to the spectral decomposition, as well as the
associated \KL{}- or proper orthogonal decomposition will be sketched.
The purpose of this short note is to show the common theoretical background
and pull some lose ends together.

%\ignore{            %%%% BEGIN IGNORE
\vspace{5mm}

\vspace{5mm}
{\noindent\textbf{MSC Classification:} \thesubject}
%}            %%%% END IGNORE

\end{abstract}

%  $Log: abstract.tex,v $
%  Revision 1.2  2018/07/05 22:01:22  hgm
%  final for arXiv
%
%  Revision 1.1  2018/06/21 23:36:24  hgm
%  first check-in old file from Lin-Param
%
%
%
%
%

%%% Local Variables: 
%%% mode: latex
%%% TeX-master: "../18_RV-algebra-model"
%%% End: 

\setcounter{page}{0}
\thispagestyle{empty}

%\pagenumbering{roman}

% table of contents if desired
%\cleardoublepage

%\tableofcontents
%\cleardoublepage

% list of figures if desired
%\listoffigures
%\cleardoublepage

%list of symbols if desired
%\listofsymbols
\cleardoublepage

\pagenumbering{arabic}

% use include to start each section on a new page.

% !TEX root = ../18_RV-algebra-model.tex
% !TEX encoding = UTF-8 Unicode
% RCSID:       $Id: introduction.tex,v 1.2.1.1 2018/09/01 01:15:28 hgm Exp $
% Author:      $Author: hgm $
% Contact:     wire@tu-bs.de
% =================================
%% texfile{
%%  AUTHOR    = "$Author: hgm $",
%%  VERSION   = "$Revision: 1.2.1.1 $",
%%  DATE      = "$Date: 2018/09/01 01:15:28 $",
%%  FILENAME  = "$RCSfile: introduction.tex,v $"}
%
% =================================

\section{Introduction}  \label{S:intro}
Probabilistic and parametric models, used in many areas of science, engineering,
and economics, share many similarities.  Probabilistic models are used to describe
uncertainties or random phenomena, whereas parametric models describe variations
or changes of some system as some parameters are changed.
Typically these are part of some larger mathematical
model describing some system with such characteristics.  A parameter can of course
be a \emph{random variable}, and this is the connection between these two
kinds of models.  Here the interest is mainly in system models with an infinite
dimensional state space, e.g.\ systems described by ordinary or partial
differential equations.  This often also makes it necessary to theoretically
consider infinitely many parameters.  In an actual numerical computation
this has of course to be reduced through some kind of discretisation to
a finite number.  And obviously one would like to have this number as small
as possible while still retaining acceptable accuracy.  This is the
realm of \emph{reduced order} models.

These reduced order models lessen the possibly
high computational demand, and are hence probabilistic or parametrised
reduced order models.
The survey \citep{BennWilcox-paramROM2015} and
the recent collection \citep{MoRePaS2015}, as well as
the references therein, provide a good account of parametric
reduced order models and some of the areas where they appear.
The interested reader may find there further information on
parametrised reduced order models and how to generate them.

Here we build on our recent work \citep{hgmRO-1-2018, hgmRO-2-2018}
analysing parametrised reduced order systems, which itself
is a continuation of \citep{boulder:2011}.  In these publications
the theoretical background of such parametrised models is treated
in a functional analysis setting.  The purpose of the present note
is to use the same kind of techniques for stochastic or probabilistic
models, where some generalisations are required due to the wish to cover
infinite dimensional state spaces, and combine this with the description
of parametric reduced order models.

As an example, assume that some physical system is investigated, which
is modelled by an evolution equation for its state $v(t) \in \C{V}$
at time $t \in [0,T]$,
where  $\C{V}$ is assumed to be a Hilbert space for the sake of
simplicity: $\dot{v}(t) = A(\vsigma, \mu;v(t)) + f(\vsigma, \mu;t);\quad v(0) = v_0$,
where the superimposed dot signifies the time derivative,
$A$ is an operator modelling the
physics of the system, and $f$ is some external excitation.
Here $\vsigma$ is a random variable (RV) defined on an event space
$\Omega$ with values in some Hilbert space $\E{S}$
(again for simplicity),
and $\mu \in \C{M}$ are parameters that can be controlled, and can be
used to evaluate the \emph{design} of the system, \emph{control}
its behaviour, or \emph{optimise} the performance in some way.
No specific structure is assumed for the set $\C{M}$.
We assume that for all possible values of $\vsigma$ and
for all $\mu$ of interest the system is well-posed.
This will make the system state $v(\vsigma,\mu;t)$ a random variable as
well, depending on the value of the parameters $\mu$.

One may be interested in the state of the system $v(\vsigma,\mu;t)$
and its statistics, or some functional of it, say $\Psi(\mu)=\D{E}(\psi(v(\vsigma,\mu))$,
where $\D{E}$ is an expectation operator.
While evaluating $A(\vsigma,\mu)$ or
$f(\vsigma,\mu)$ for a certain $\mu$ may be straightforward,
evaluating $v(\vsigma,\mu;t)$ or $\Psi(\mu)$ may be very costly.
This is why one wants representations of $v(\vsigma,\mu;t)$ or $\Psi(\mu)$
which allow a cheaper evaluation.  This is achieved through reduced order models,
which are often also called \emph{proxy}- or \emph{surrogate}-models.
It turns out that such random and parametric objects
can be analysed by associated linear maps \citep{hgmRO-1-2018, hgmRO-2-2018},
which renders them much more accessible to the techniques of linear functional
analysis, a well understood subject.
This association with linear mappings has probably been
known for a long time, see \citep{kreeSoize86} for an exposition
in the context of stochastic models.  In \refS{parametric} the association
of parametric and stochastic models with linear maps will be explained,
in passing touching on reproducing kernel Hilbert spaces.  
The classical probabilistic framework (cf.\ \citep{Sullivan2015}), starting
from measurable spaces and $\sigma$-algebras, can be used to define 
algebras of random variables (RVs) as measurable functions on these measure
spaces, and the expectation operator as integral of these RVs w.r.t.\ the
probability measure.  These algebras of RVs can be used in the case of
probabilistic models to build the range or image space for these linear maps
as spaces of classical RVs.  But alternatively one may also start by using
as fundamental concepts algebras of objects that we want to call RVs
together with the expectation operator (cf.\ \citep{Segal1978}) as a linear
functional, and if this algebra of RVs is Abelian or commutative one essentially
recovers equivalence with classical probability.  This approach allows for
non-commuting algebras of RVs, which is important (cf.\ \citep{MingoSpeicher2017})
in order to deal with e.g.\ random matrices, random fields of tensors,
quantum theory and quantum fields.  More important for our
immediate purposes here, this view greatly facilitates the specification of
stochastic models on infinite dimensional spaces.  Such an algebra of RVs,
whether generated classically as derived concept as an algebra of measurable
functions, or used as a primary model of possibly non-commuting of RVs, seems to be a
natural object to use in the case of stochastic models on infinite dimensional vector
spaces, as it allows to generalise such stochastic models to so-called
weak distributions or generalised processes (cf.\ \citep{segal58-TAMS, LGross1962, 
gelfand64-vol4, segalNonlin1969}), and thereby elegantly circumvent many
problems which arise when one tries to define $\sigma$-additive set functions
for example on Hilbert spaces.
This algebraic and analytic view on probability will be explained in \refS{alg-RV}.
Everything is tied together in \refS{correlat} in the analysis of the
generalised correlation operator, its factorisations, as well as its
spectral decomposition, and the last \refS{concl} concludes by pointing out
once more the connection between functions in high-dimensional spaces and
the associated linear maps and correlation operators, where well-known methods
can be used to analyse their structure.

%  $Log: introduction.tex,v $
%  Revision 1.2.1.1  2018/09/01 01:15:28  hgm
%  after revision
%
%  Revision 1.2  2018/07/05 22:01:30  hgm
%  final for arXiv
%
%  Revision 1.1  2018/06/21 23:37:47  hgm
%  first check-in old file from Lin-Param
%
%
%
%
%
%

%%% Local Variables: 
%%% mode: latex
%%% TeX-master: "../18_RV-algebra-model"
%%% End: 

% !TEX root = ../18_RV-algebra-model.tex
% !TEX encoding = UTF-8 Unicode
% RCSID:       $Id: param-map.tex,v 1.2.1.1 2018/09/01 01:17:41 hgm Exp $
% Author:      Hermann G. Matthies
% Contact:     wire@tu-bs.de
% =================================

\section{Parametric and stochastic models} \label{S:parametric}
We start with a short recap of \citep{hgmRO-1-2018, hgmRO-2-2018}, where
the interested reader may find more detail.  Let $r: \C{M} \to \C{U}$ be
a generic substitute for any
one of the parametric objects alluded to in the introduction, e.g.\ things like
$\mu \mapsto v(\vsigma,\mu,t)\in\C{V}$ or $\mu \mapsto \dot{v}(\vsigma,\mu,\cdot)\in
\Lp_2([0,T])\otimes\C{V}$; $\omega \mapsto \vsigma(\omega)\in\E{S}$---with $\Omega$
taking the rôle of $\C{M}$; $(\mu,\omega) \mapsto v(\vsigma(\omega),\mu,t)\in\C{V}$
---with $\C{M}\times\Omega$ taking the rôle of $\C{M}$; 
$\omega \mapsto f(\vsigma(\omega),\mu,t)\in\C{V}$---with $\Omega$ taking the rôle of $\C{M}$,
or $\mu \mapsto A(\vsigma,\mu,\cdot)\in(\C{V}\to\C{V})$---the space of maps from
$\C{V}$ to $\C{V}$, etc.
 
The space $\C{U}$ is assumed for the sake of simplicity
as a separable Hilbert space.  The function $r$ can thus be either
a parametric input, or a random input---i.e.\ a random variable (RV), in
which case $\C{M}$ would be a measure space---to a model like that described
in \refS{intro}, or the operator of that model, or the state (solution)
of that system.
Assuming---without significant loss of generality---that the image
$\spn r(\C{M}) = \spn \im r \subseteq \C{U} $ is dense in $\C{U}$,
one may to each such function $r$ associate a linear map
$R: \C{U} \ni u \mapsto \bkt{r(\cdot)}{u}_{\C{U}} \in \D{R}^{\C{M}}$ into
the space $(\C{M}\to\D{R})$ of all real-valued functions on $\C{M}$.
By construction, $R$ restricted to $\spn\im r = \spn r(\C{M})$ is injective.
In \refS{alg-RV} it will be explained how---in the case of a probabilistic or
random model---the Hilbert space can be generated from an algebra of RVs.

As an aside, note that
on its restricted range $\tilde{\C{R}}:=R(\spn \im r) \subseteq \D{R}^{\C{M}}$
one may define an inner product as $\bkt{\phi}{\psi}_{\C{R}} :=
\bkt{R^{-1} \phi}{R^{-1} \psi}_{\C{U}}$ for all $\phi, \psi \in \tilde{\C{R}}$.
Denote the completion with this inner product by $\C{R}$.  This makes
$R$ and $R^{-1}$ into bijective isometries, hence \emph{unitary} maps
between $\C{U}$ and $\C{R}$.  It may easily be shown \citep{hgmRO-1-2018, hgmRO-2-2018}
that $\C{R}$ is a \emph{reproducing kernel Hilbert space} (RKHS)
\citep{berlinet, Janson1997} with reproducing kernel
$\vkappa(\mu_1, \mu_2) := \bkt{r(\mu_1)}{r(\mu_2)}_{\C{U}}$,
such that the reproducing property
$\bkt{\vkappa(\mu,\cdot)}{\phi}_{\C{R}} = \phi(\mu)$ holds
for all $\phi \in \C{R}$.
In this note the RKHS $\C{R}$ will not be used, but the 
important thing to keep in mind is that the map $R$ and the space $\C{R}$
of scalar functions on the set $\C{M}$---one might view them as problem
oriented co-ordinates---carry the same information as the parametric
object $r(\mu)$.

Often some information of what is important in the set $\C{M}$ is also
available, here it is assumed to be given by a Hilbert subspace
$\C{Q}\subseteq \D{R}^{\C{M}}$, usually different from $\C{R}$.
From now on we shall by slight abuse of notation view the map $R$
as mapping into $\C{Q}$ and still assume that it is injective
as well as closed, for the sake of simplicity.
Details like the assumption that the subspace
$R^{-1}(\C{Q})$ is dense in $\C{U}$ will not
always be spelt out in detail for the sake of brevity.
The idea is that with $u\in\C{U}$ of unit length the vectors  $R u\in \C{Q}$
with large norm are more important, and this will be considered in
building reduced order models.  As will be shown \citep{hgmRO-1-2018, hgmRO-2-2018}
in \refS{correlat}, the map $C:\C{U}\to\C{U}$ defined by $C=R^* R$, where $R^*$
is the adjoint of $R$, is central to the analysis.  More precisely,
with the above assumptions on $R$ the adjoint $R^*$ is surjective,
and $C$ is a densely defined self-adjoint positive definite operator,
which we shall call the `correlation' of the model $r(\mu)$. 

A random variable or stochastic model as exemplified by
the RV $\vsigma$ in \refS{intro} is usually formulated
as a measurable map $\vsigma:\Omega\to\E{S}$, where $(\Omega, \F{A}, \D{P})$
is a probability space with $\sigma$-algebra $\F{A}$ and
probability measure $\D{P}$.  One may view the set $\Omega$ as a parameter
set like $\C{M}$ above, and one can construct a linear map into
the space $\D{R}^\Omega$, i.e.\ the scalar random variables.
Without loss of generality, we assume that
$\spn \vsigma(\Omega) = \spn \im \vsigma \subseteq \E{S}$ is dense in
the separable Hilbert space $\E{S}$, and define  \citep{kreeSoize86}
\begin{equation}   \label{eq:lin-m-st}
 S: \E{S} \ni \xi \mapsto \bkt{\vsigma(\cdot)}{\xi}_{\E{S}} \in \D{R}^{\Omega} .
\end{equation}
It remains to define an inner product on $\D{R}^\Omega$ and a subspace
corresponding to $\C{Q}$ for the parametric case above.  This will be done in
\refS{alg-RV}.  For the time being assume that this has been defined, i.e.\ there
is an inner product $\bkt{\cdot}{\cdot}_{\E{V}}$ and a corresponding
Hilbert space of (equivalence classes) of RVs $\E{V}\subseteq \D{R}^{\Omega}$,
and we regard $S$ as a map $S:\E{S}\to\E{V}$ with the same properties
as assumed for $R$ above.  Obviously the densely
defined self-adjoint positive definite operator $C_{\vsigma} = S^* S : \E{S}\to\E{S}$
corresponding to $C=R^* R$ above is indeed the correlation operator
of the RV $\vsigma$.

In case $\vsigma$ is an input to a dynamical system like the one alluded
to in \refS{intro}, the state of the system $v(\vsigma,\mu;t)$ also becomes
a stochastic quantity, and inner product with a vector $w\in\C{V}$
leads for fixed $\mu$ and $t$ automatically to a linear mapping
\begin{equation}   \label{eq:lin-u-st}
 P: \C{V} \ni w \mapsto \bkt{v(\vsigma(\cdot),\mu;t)}{w}_{\C{V}} \in \D{R}^{\Omega} ,
\end{equation}
which we shall regard again as a map $P:\C{V}\to\E{V}$ into the just defined
space $\E{V}$.  This defines a third correlation operator
$C_v = P^* P:\C{V}\to\C{V}$.

It may be seen that with the correspondences
\begin{equation}   \label{eq:lin-correspond}
R:\C{U} \to \C{Q} \quad \text{---} \quad S:\E{S} \to \E{V} \quad \text{---} \quad
P:\C{V}\to\E{V}
\end{equation}
all three situations are completely analogous, and may in the simplest case be
dealt with in the same formalism.
The idea on how to obtain representations of $r(\mu)$ resp.\ $\vsigma(\omega)$
resp.\ $v(\omega)$
is the following \citep{hgmRO-1-2018, hgmRO-2-2018}, which we shall mainly
demonstrate for $r(\mu)$:  choose a complete basis
$\{q_j \}_j \subset \C{Q}$, and represent $r(\mu)$ as
\begin{equation}   \label{eq:basic-repr}
r(\mu) = \sum_j \alpha_j R^* q_j(\mu) .
\end{equation}
A good reduced order model is one where
\begin{equation}   \label{eq:basic-repr-ROM}
 r_{\text{ROM}}(\mu)  = \sum_{j=1}^J \alpha_j R^* q_j(\mu) 
\end{equation}
is a good approximation to $r(\mu)\approx r_{\text{ROM}}(\mu)$ with
a small $J$, i.e.\ with not too many terms.  In \refS{correlat}
some other possibilities for the choice of basis $\{q_j \}_j$ will be discussed,
where the $\mu$-dependence is encoded in the scalar functions from $\C{Q}$,
but where a basis of $\mu$-independent vectors is picked from $\C{U}$, and
where again for the sake of brevity and simplicity we shall confine ourselves
to complete orthonormal systems (CONS).  The important message here is that
with $R$ one has a factorisation of $C=R^* R$, and that the adjoint is
the map which carries a representation on the function space to
the space $\C{U}$.  Later we shall indicate \citep{hgmRO-1-2018, hgmRO-2-2018}
how every representation leads to a factorisation of $C$, and that
--- with some additional assumptions on $C$ --- every factorisation leads
to a representation.  But the description and analysis via factorisations
is more general \citep{segal58-TAMS, LGross1962, gelfand64-vol4, segalNonlin1969,
kreeSoize86}, and this is needed in the formulation of probabilistic models
where $\C{U}$ resp.\ $\E{S}$ is an infinite dimensional Hilbert space.

%  $Log: param-map.tex,v $
%  Revision 1.2.1.1  2018/09/01 01:17:41  hgm
%  after revision
%
%  Revision 1.2  2018/07/05 22:01:52  hgm
%  final for arXiv
%
%  Revision 1.1  2018/06/21 23:38:30  hgm
%  first check-in old file from Lin-Param
%
%
%
%

%%% Local Variables: 
%%% mode: latex
%%% TeX-master: "../18_RV-algebra-model"
%%% End: 

% !TEX root = ../18_RV-algebra-model.tex
% !TEX encoding = UTF-8 Unicode
% RCSID:       $Id: rv-algebra.tex,v 1.2.1.3 2018/09/01 22:26:30 hgm Exp $
% Author:      $Author: hgm $
% Contact:     wire@tu-bs.de
% =================================
%% texfile{
%%  AUTHOR    = "$Author: hgm $",
%%  VERSION   = "$Revision: 1.2.1.3 $",
%%  DATE      = "$Date: 2018/09/01 22:26:30 $",
%%  FILENAME  = "$RCSfile: rv-algebra.tex,v $"}
%
% =================================

\section{Algebras of random variables}  \label{S:alg-RV}
Here we shall take a closer look at the stochastic or probabilistic model
$\vsigma:\Omega\to\E{S}$ and the associated linear map $S:\E{S}\to\E{V}$,
as well as the space of RVs $\E{V}$ and how it is generated.
Although there are classical ways of specifying the space $\E{V}$,
the most natural one seems to be the algebraic approach to probability.
These ideas are certainly also used in the classical approach,
but the algebraic probability approach distills the essential components
in an abstract setting and allows at the same time generalisations.
Historically, when looking back as how in the beginnings of probability
theory the Bernoullis treated random variables (RVs), it is clear that
they added them and took multiples---hence they form a vector space---and
that they multiplied them with each other---so they form an algebra.
Although the formalisation of probability as formulated by Kolmogorov
used the concept of measure and this algebraic background was
largely ignored, it was revived with the advent of quantum theory.
It turns out that here this view is essential, as not all observables
can be observed simultaneously, and this is reflected in the fact
that they do not commute in the algebra.  Another topic where this
view is very advantageous are random matrices and more generally
random fields of even-order tensors.

We are mainly interested in `real' or self-adjoint RVs as they will later be
called.  But for analytical convenience we shall treat complex RVs,
following Paul Painlevé's and Jacques Hadamard's adage that the shortest
path between two truths in the real domain passes through the complex domain
--- ``le plus court chemin entre deux vérités dans le domaine
réel passe par le domaine complexe''.
Some algebraic language is needed, but
most of the terms will be familiar from complex numbers and from matrices,
which are indeed two simple but prime examples of algebras.
Let us start right away with a simple and mostly familiar example
from probability theory, which will at the same time serve as
motivation, concrete example, and explanation of the abstract setting.

\subsection{Specifying the algebra}  \label{SS:specify-alg}
Consider a probability space $(\Omega,\F{A},\D{P})$ with a set of elementary
events $\Omega$, $\sigma$-algebra $\F{A}$ of measurable subset of $\Omega$,
and probability measure $\D{P}$.
In the vector space $\Lp_0(\Omega,\F{A},\D{P};\D{C})$ of complex-valued
measurable functions / classical random variables on $\Omega$---which for the
sake of brevity shall be denoted just by $\Lp_0(\Omega)$---let
$\C{A}_s:=\Lp_{0s}(\Omega)\subset\Lp_0(\Omega)$
be the vector subspace of complex-valued simple measurable functions,
i.e.\ complex linear combinations of functions
$\mat{1}_{\C{E}}$, which for $\C{E}\in\F{A}$ are defined to be
$\mat{1}_{\C{E}}(\omega)=1$ if $\omega\in\C{E}\subseteq\Omega$, and
zero otherwise.  Hence $\C{A}_s$ are the RVs where each one of them
can only take finitely many different values.

On this vector space we may define a multiplication by just pointwise
multiplication of two such RVs, and the product is obviously again a
simple function; in fact for $\C{E}, \C{F}\in\F{A}$ one has
$\mat{1}_{\C{E}} \mat{1}_{\C{F}} = \mat{1}_{\C{E}\cap\C{F}}$, i.e.\ the
multiplication in $\C{A}_s$ reflects the intersection in the $\sigma$-algebra
$\F{A}$.  This means that the space $\C{A}_s$ is closed under
multiplication and hence thanks to the properties of the multiplication
on $\D{C}$ is a complex, associative, and commutative or Abelian algebra,
with the familiar distributive law from $\D{C}$ coupling addition and multiplication
also on $\C{A}_s$.  Another way of saying this is to state that the multiplication
is a bilinear map from $\C{A}_s\times\C{A}_s$ to $\C{A}_s$.
Let us note in passing that with the same definition of pointwise multiplication
also $\Lp_0(\Omega)$ is an associative and commutative algebra---with $\C{A}_s$
a sub-algebra---as the pointwise product of two measurable functions is again
measurable, but we shall see later that for our purposes $\Lp_0(\Omega)$ is in
general too big.  The element $\mat{1}_{\Omega}\in\C{A}_s\subset\Lp_0(\Omega)$
which is constant equal to unity is obviously a neutral element or unit for
the multiplication, and hence $\C{A}_s$ and $\Lp_0(\Omega)$ are called unital algebras.
For $\psi\in\C{A}_s$ one can now compute \emph{powers} $\psi^n = \psi \psi^{n-1}$ for any
integer $n\ge 1$, and if we define $\psi^0=\mat{1}_{\Omega}$ in a unital
algebra even for any $n\ge 0$.  Given a polynomial $Q(X)=\sum_{k=0}^n \alpha_k X^k\in\Pi_1$
in one unknown $X$ with complex co-efficients $\alpha_k\in\D{C}$,
it is now possible to evaluate $Q(\psi)\in\C{A}_s$ for any $\psi\in\C{A}_s$.
For some $\phi\in\C{A}_s$ there is a $\psi\in\C{A}_s$ such that $\phi\psi=\mat{1}_{\Omega}$.
This is then called the (multiplicative) \emph{inverse} $\psi=\phi^{-1}$,
 such that $\phi \phi^{-1} = \mat{1}_{\Omega}$.

For a complex number $\zeta\in\D{C}$ its complex conjugate is denoted by
$\zeta^*\in\D{C}$, and this operation is an \emph{involution}, as $(\zeta^*)^*=\zeta$.
One may extend this involution from $\D{C}$ to the algebra $\Lp_0(\Omega)$ through
a pointwise definition of complex conjugation, and hence also to its sub-algebra $\C{A}_s$.
For $\phi, \psi\in\Lp_0(\Omega)$ and $\zeta\in\D{C}$ this involution obviously
satisfies $(\phi + \zeta\psi)^* = \phi^* + \zeta^* \psi^*$ and is thus \emph{anti-linear}.
As regards the product of two RVs, it satisfies $(\phi\psi)^* = \psi^*\phi^*$,
and it is easy to verify that both $\C{A}_s$ and $\Lp_0(\Omega)$ are closed under this involution.
Associative algebras with such an anti-linear involution and the indicated behaviour on
products are called $^*$-algebras---the element $\psi^*$
is usually called in algebraic terms the \emph{adjoint} of $\psi$---and both
$\Lp_0(\Omega)$ and its sub-algebra $\C{A}_s=\Lp_{0s}(\Omega)$ are thus $^*$-algebras.

Let $\Pi_2^c$ denote the set of all polynomials $Q(X,Y)$ with complex co-efficients
in two commuting variables $X, Y$.  For $\phi\in\C{A}_s$
the unital sub-$^*$-algebra $\D{C}[\phi,\phi^*]:=\{Q(\phi,\phi^*)\mid Q\in \Pi_2^c \} 
\subset\C{A}_s$ is called the sub-algebra \emph{generated} by $\phi\in\C{A}_s$.
Observe that if $\psi\in\Lp_0(\Omega)$ is \emph{self-adjoint}, i.e.\ $\psi=\psi^*$,
then $\psi$ has only real values, and if $\psi = \phi^*\phi$ for some
$\phi\in\Lp_0(\Omega)$, then $\psi$ is self-adjoint (real) and is called
\emph{positive} as it can not take negative values,
i.e.\ $0\le \psi = \phi^*\phi$---in case $0<\psi$ it is usually
called \emph{strictly positive}.
One says that for self-adjoint $\phi, \psi\in\Lp_0(\Omega)$ one has
$\psi \le \phi$ iff $\phi-\psi$ is positive, and thus one can define a
partial order on $\C{A}_s$ and $\Lp_0(\Omega)$.
Positive self-adjoint elements $\psi\in\Lp_0(\Omega)$ which are \emph{idempotent},
i.e.\ satisfy $\psi^2 = \psi \psi = \psi$, are called \emph{projections}.  Observe that
each $\mat{1}_{\C{E}}$ is a projection, and that the unit $\mat{1}_{\Omega}$ is a
maximal projection in the order mentioned.
In fact all projections in $\Lp_0(\Omega)$ and $\C{A}_s$
have the form $\mat{1}_{\C{E}}$ for some $\C{E}\in\F{A}$.
Ultimately, one is only interested in the self-adjoint elements of the algebra
$\C{A}_s$, as they take real values; they are therefore often also called \emph{observables}.
The other elements of the algebra may be regarded as merely a kind of
analytical completion to make the theory nice.  It may be remarked that
the self-adjoint elements of $\C{A}_s$ form a \emph{real} subspace of $\C{A}_s$.
Obviously an arbitrary $\phi\in\C{A}_s$ may be decomposed into real
and imaginary parts: $\phi = \Re \phi + \ii \Im \phi$ with real resp.\ self-adjoint
$\Re \phi = {(\phi + \phi^*)}/{2}$ and $\Im \phi = (\phi - \phi^*)/(2 \ii)$,
so that the whole algebra is the complex span of the self-adjoint elements
or observables.

To extract the essential point from this example and generalise,
we start with an associative algebra $\C{A}$ of what we want to call random
variables (RVs) $\tns{a}, \tns{b},\dots \in \C{A}$, i.e.\ a vector space
\citep{Segal1978} equipped with an associative and bi-linear multiplication
which will be denoted just by juxtaposition: $\C{A}\times\C{A} \ni (\tns{a},\tns{b})
\mapsto \tns{ab} \in \C{A}$.  As was noted before, it is advantageous to assume
the algebra to be a complex algebra, which
is no loss of generality as any real algebra may be embedded into a
complex one.  For $\tns{a}\in\C{A}$ the powers $\tns{a}^n$ are defined for any
integer $n\ge 1$ in the natural recursive fashion.
Additionally assume that the algebra is unital, i.e.\ has a multiplicative
unit $\tns{e}$ such that $\tns{ae}=\tns{ea}=\tns{a}$ for any $\tns{a}\in\C{A}$,
and one defines the power $\tns{a}^n$ for $n=0$ by $\tns{a}^0=\tns{e}$.
Hence for a polynomial $Q(X)\in\Pi_1$ it is now possible
to evaluate $Q(\tns{a})$ for any $\tns{a}\in\C{A}$.  Also assume
that there is an anti-linear involution defined,
called the `adjoint', denoted as $\tns{a}^*$,
such that $(\tns{a}^*)^* = \tns{a}$ and $(\tns{ab})^* = \tns{b}^* \tns{a}^*$.

Let $\Pi_2^n$ be a set of all polynomials $Q(\tns{X},\tns{Y})$ with complex
co-efficients in two non-commuting variables $\tns{X},\tns{Y}$, then for $\tns{a}\in\C{A}$
the unital sub-$^*$-algebra
$\D{C}\{\tns{a},\tns{a}^*\}:=\{Q(\tns{a},\tns{a}^*)\mid Q\in \Pi_2^n \}
\subset\C{A}$ is called the sub-algebra generated by $\tns{a}\in\C{A}$.  Elements
$\tns{a}\in\C{A}$ such that $\tns{a}=\tns{a}^*$ are called self-adjoint,
and self-adjoint elements which may be factored as $\tns{a}=\tns{b}^*\tns{b}$
are called positive.  Positive elements form a salient pointed cone
which defines an order relation on $\C{A}$.  Positive elements $\tns{p}$ which are
idempotent $\tns{p}=\tns{p}\tns{p}=\tns{p}^2=\tns{p}^*\tns{p}=\tns{p}^*$ are called projections.
Observe that $\tns{e}$ is a projection, and that it is maximal w.r.t.\ the order
mentioned.  Succinctly stated, we assume that $\C{A}$ is a complex
associative unital $^*$-algebra, not \emph{necessarily} commutative.
As was shown, both $\Lp_0(\Omega)$ and $\C{A}_s$ considered above are
commutative examples of such algebras.  Again, one is later ultimately
interested in the self-adjoint elements of $\C{A}$---the \emph{observables}.
Also in the general abstract case they form a \emph{real} subspace of $\C{A}$,
and an arbitrary $\tns{a}\in\C{A}$ may be decomposed into two parts
$\tns{a} = \tns{a}_s + \ii \tns{a}_w$ with self-adjoint $\tns{a}_s = 
(\tns{a} + \tns{a}^*)/2$ and $\tns{a}_w = (\tns{a} - \tns{a}^*)/(2 \ii)$
---also called the \emph{symmetric} and \emph{skew} parts---so that the
whole algebra is the complex span of the self-adjoint elements, the observables.
And naturally, if for some $\tns{a}\in\C{A}$ there is a $\tns{c}\in\C{A}$ such that
$\tns{ac} = \tns{ca} = \tns{e}$, then $\tns{c}=\tns{a}^{-1}$
is the unique multiplicative inverse of $\tns{a}$.

\subsection{States and the expectation functional}   \label{SS:exp-func}
To continue, we return to the example $\C{A}_s$ above.  Just as classical
probability builds on the measurable space $(\Omega,\F{A})$ on one hand
and the probability measure $\D{P}$ on the other hand, in the algebraic framework
the second entity needed is the linear expectation functional $\D{E}:\C{A}_s\to\D{C}$.
To define the expected value for a RV $\phi\in\C{A}_s$ one only
has to look at the generating elements $\mat{1}_{\C{E}}$ with $\C{E}\in\F{A}$.
Here one defines $\EXP{\mat{1}_{\C{E}}} := \int_\Omega \mat{1}_{\C{E}}(\omega)\,
\D{P}(\di \omega) = \D{P}(\C{E})$ and extends this
by linearity to all of $\C{A}_s$.  Thus the probability of an event
$\C{E}\in\F{A}$ is given in terms of the expected value of the associated
projection $\mat{1}_{\C{E}}$.  For a typical 
$\phi(\omega) = \sum_k \alpha_k \mat{1}_{\C{E}_k}(\omega) \in  \C{A}_s$
with $\alpha_k\in\D{C}$ this gives $\EXP{\phi} = \int_{\Omega} \phi(\omega)\, 
\D{P}(\di \omega) = \sum_k \alpha_k \D{P}(\C{E}_k)\in\D{C}$.
Obviously, as $\D{P}(\Omega)=1$, the expected value of the unit is
$\EXP{\mat{1}_{\Omega}} = 1$, a kind of \emph{normalisation} of the
expectation functional.

This linear functional $\D{E}$ additionally satisfies
$\EXP{\phi^*} = (\EXP{\phi})^*$ and thus carries the adjoint to its
complex conjugate and hence is real on self-adjoint elements.  Such
a linear functional is itself called \emph{self-adjoint}.  In addition,
$\EXP{\phi^* \phi} = \sum_k (\alpha_k^*\alpha_k) \EXP{\mat{1}_{\C{E}_k}}
= \sum_k |\alpha_k|^2 \D{P}(\C{E}_k)\ge 0$, i.e.\  the functional
is non-negative on positive $\psi=\phi^* \phi \in \C{A}_s$.
Such a self-adjoint linear functional is itself called \emph{positive}.
If $\rho\in \C{A}_s$ is positive with unit expected value $\EXP{\rho}=1$,
one may define a new expectation functional---corresponding to a change of
probability measure---via $\D{E}_{\rho}(\phi) := \EXP{\rho \phi}= 
\int_\Omega \rho(\omega)\,\phi(\omega)\, \D{P}(\di \omega) $.
It is easily checked that $\D{E}_{\rho}$ is linear, self-adjoint,
positive, and normalised.  Such linear functionals which can serve as
expectation are called \emph{states}, an element of the \emph{dual}
space $\C{A}_s^*$.

The element $\bar{\phi} := \EXP{\phi} \mat{1}_{\Omega}
\in \C{A}_s$ is called the \emph{mean} of $\phi\in\C{A}_s$ and the additive rest
$\tilde{\phi} = \phi - \bar{\phi}\in\C{A}_s$ is its \emph{zero-mean} or \emph{centred}
or \emph{fluctuating} part.  The one-dimensional unital $^*$-algebra
$\C{A}_{sc}:=\D{C}[\mat{1}_{\Omega}]=\spn\{\mat{1}_{\Omega}\} \subset\C{A}_s$---isomorphic
to $\D{C}$---are the \emph{constants}, whereas the subspace
$\C{A}_{s0} := \ker \D{E}$ are the \emph{zero-mean} or \emph{centred} RVs, such that
$\C{A}_s = \C{A}_{sc} \oplus \C{A}_{s0}= \D{C}[\mat{1}_{\Omega}]\oplus\ker\D{E}$
as a direct sum.

One may observe that in general not every
measurable $\phi\in\Lp_0(\Omega)$ has a finite integral.  
Thus the algebra of \emph{all} classical RVs $\Lp_0(\Omega)$
is too big for our purpose as one would like $\EXP{\cdot}$
to be defined on the whole algebra.  This is the reason to start with the
`smaller' algebra $\C{A}_s=\Lp_{0s}(\Omega)$.  It is a building block from
which more complicated RVs can be built via limiting processes.

In the general abstract case one also wants a linear, self-adjoint, positive,
and normalised functional---a \emph{state}---$\D{E}:\C{A}\to\D{C}$ with
$\EXP{\tns{a}^*} = \EXP{\tns{a}}^*$.  Such a state is called \emph{faithful}
if $\EXP{\tns{a}^* \tns{a}}=0$ implies $\tns{a}=0$.
If a state is not faithful, then one can start to work with an algebra
of equivalence classes, where two elements $\tns{a}, \tns{b} \in \C{A}$ are considered
equivalent iff $\EXP{(\tns{a}-\tns{b})^*(\tns{a}-\tns{b})}=0$.  It is therefore no loss of
generality to assume that the state is faithful.
The projections $\tns{p}\in\C{A}$ are
also identified with events, and the \emph{probability} of the
event $\tns{p}\in\C{A}$ may be \emph{defined} as $\D{P}(\tns{p}):=\EXP{\tns{p}}$.
As $\D{E}$ is positive, one has $\D{P}(\tns{p})\ge 0$, and as
$\tns{e}$ is a maximal projection, $\D{P}(\tns{p}) \le \D{P}(\tns{e})=1$.
One defines the mean part of a RV
as a multiple of the identity $\bar{\tns{a}}:= \EXP{\tns{a}} \tns{e}$ and
the fluctuating zero-mean or centred part as $\tilde{\tns{a}}:= \tns{a} - 
\bar{\tns{a}}$ with $\EXP{\tilde{\tns{a}}}=0$.
The one dimensional sub-$^*$-algebra $\C{A}_c=\D{C}[\tns{e}]=\spn\{\tns{e}\}$ of
constants---isomorphic to $\D{C}$---are multiples of the identity,
and the subspace of zero-mean fluctuating parts
$\C{A}_0=\ker \D{E}$ is the kernel of the state, and the whole algebra
is the direct sum of both parts $\C{A} = \C{A}_c \oplus \C{A}_0 =
\D{C}[\tns{e}]\oplus\ker \D{E}$.  An abstract algebra which satisfies
all these requirements together with a distinguished faithful state
as expectation is called a
\emph{probability algebra}.  If $\tns{\vrho}\in\C{A}$ is positive
with unit expectation $\EXP{\tns{\vrho}}=1$, then one may define
a new weighted state by $\D{E}_{\tns{\vrho}}(\tns{a}):=\EXP{\tns{\vrho}\,\tns{a}}$
for $\tns{a}\in\C{A}$.

A faithful state may be used to define an inner product on $\C{A}$
\citep{Segal1978, segal58-TAMS, segalNonlin1969} via a positive definite
sesqui-linear form:
\begin{equation}  \label{eq:XXIII-inner-p}
  \C{A}^2 \ni (\tns{a},\tns{b}) \mapsto \bkt{\tns{a}}{\tns{b}}_{2} :=
  \EXP{\tns{b}^* \tns{a}} \in \D{C} .
\end{equation}
As usual, one may
define the square of a norm via $\nd{\tns{a}}_{2}^2:= \bkt{\tns{a}}{\tns{a}}_{2}$.
The completion of $\C{A}$ in the uniform topology generated by this norm is a
Hilbert space denoted by $\Lp_2(\C{A})$, which is one candidate for $\E{V}:=\Lp_2(\C{A})$.
Later we shall see more possible ways of generating a Hilbert space of RVs.
With this inner product the above direct sum $\C{A} = \C{A}_c
\oplus \C{A}_0$ is an \emph{orthogonal} direct sum, i.e.\ $\C{A}_c=
\D{C}[\tns{e}]=\spn \{ \tns{e}\} = (\ker \D{E})^\perp = \C{A}_0^\perp$.

As the expectation or state is normally also continuous in the topology
of the associated Hilbert space $\E{V}$, it can be defined also on $\E{V}$
giving an orthogonal decomposition $\E{V}=\ker \D{E} \oplus (\ker \D{E})^\perp
=: \E{V}_0 \oplus \D{C}[\tns{e}]$.  For the probabilistic model
$S:\E{S}\to\E{V}$ this means that it can be extended to $\xi\in\E{S}$
as $\D{E}_{\E{S}}(\xi) := \D{E}(S\xi)$,
and with it an orthogonal decomposition of $\E{S} = \E{S}_0 \oplus \E{S}_0^\perp
:= \ker \D{E}_{\E{S}} \oplus (\ker \D{E}_{\E{S}})^\perp$, where
$(\ker \D{E}_{\E{S}})^\perp= \spn\{S^*\tns{e}\}$ are multiples of the \emph{mean}
$\bar{\vsigma}:= S^*\tns{e}\in\E{S}$ of the RV $\vsigma$.  Instead of looking at
the correlation operator $C_\vsigma =  S^* S$, one is usually only
interested in the correlation $\tilde{C}_\vsigma=\tilde{S}^*\tilde{S}$ of
$\tilde{S}$, where $\tilde{S}:\E{S}\ni\xi \mapsto  S\xi - \D{E}_{\E{S}}(\xi)\tns{e}
\in\E{V}_0$---$\tilde{C}_\vsigma$ is called the \emph{covariance} operator.
Completely analogous statements can be made
for the map $P:\C{V}\ni w\mapsto\bkt{v(\vsigma)}{w}_{\C{U}}\in\E{V}$,
the associated expectation $\D{E}_{\C{V}}(w) :=\D{E}(Pw)$,
the orthogonal split $\C{V}=\C{V}_0\oplus\C{V}_0^\perp := \ker \D{E}_{\C{V}}
\oplus \spn\{P^*\tns{e}\}$, and the associated covariance operator.

In the example algebra $\C{A}_s = \Lp_{0s}(\Omega)$ from above,
identifying $\mat{1}_{\C{E}}$ and $\mat{1}_{\C{F}}$ if $\C{E}, \C{F}\in\F{A}$
differ only by a null-set $\C{N}\in\F{A}$ with $\D{P}(\C{N})=0$, the integral
or expected value becomes a faithful state.  As is well known \citep{Segal1978},
the construction in \refeq{eq:XXIII-inner-p} defines the $\Lp_2$ inner product
$\bkt{\phi}{\psi}_{2} = \EXP{\psi^*\phi}=
\int_\Omega \psi(\omega)^* \phi(\omega)\, \D{P}(\di\omega)$
for $\phi, \psi\in\C{A}_s= \Lp_{0s}(\Omega)$, and the completion is the familiar
Hilbert space $\Lp_2(\Omega) = \Lp_2(\C{A}_s)$.  The inner product 
$\bkt{\phi}{\psi}_{2}$ of two RVs $\phi, \psi\in\C{A}_s$ is also called their
\emph{correlation}, and one may continue and
define the \emph{covariance} in the usual way by $\cov(\phi,\psi) :=
\EXP{\tilde{\psi}^*\tilde{\phi}} =\bkt{\tilde{\phi}}{\tilde{\psi}}_{2}$,
i.e.\ the inner product or correlation of the fluctuating parts.
The variance of a RV $\phi\in\C{A}_s$ is then $\var(\phi):=\cov(\phi,\phi)$,
and one has from Pythagoras's theorem $\nd{\phi}_2^2 = \nd{\bar{\phi}}_2^2 +
\nd{\tilde{\phi}}_2^2 =\EXP{\phi}^2 + \var(\phi)$.
Two RVs $\phi, \psi \in\C{A}_s$ are \emph{uncorrelated} iff their
covariance vanishes: $\cov(\phi,\psi) = 0$, i.e.\ their fluctuating parts
are \emph{orthogonal}.  Two such RVs are \emph{independent} iff
$\cov(Q_1(\phi,\phi^*),Q_2(\psi,\psi^*)) = 0$ for all $Q_1, Q_2 \in \Pi_2^c$ with
$\EXP{Q_1(\phi,\phi^*)}= \EXP{Q_2(\phi,\phi^*)}=0$, i.e.\ if the centred
subspaces of the algebras generated by them are orthogonal,
i.e.\ $(\D{C}[\phi,\phi^*]\cap\ker\D{E}) \perp (\D{C}[\psi,\psi^*]\cap\ker\D{E})$.

Completely analogous in the general case, for two RVs $\tns{a}, \tns{b} \in\C{A}$
one defines the \emph{correlation} as the inner product $\bkt{\tns{a}}{\tns{b}}_{2}$,
the \emph{covariance} as the inner product of the fluctuating parts
$\cov(\tns{a},\tns{b}) := \bkt{\tilde{\tns{a}}}{\tilde{\tns{b}}}_{2}$, and the
\emph{variance} as $\var(\tns{a}):=\cov(\tns{a},\tns{a})$.  
Pythagoras's theorem can be applied here as well to give $\nd{\tns{a}}_2^2 = 
\nd{\bar{\tns{a}}}_2^2 + \nd{\tilde{\tns{a}}}_2^2 =\EXP{\tns{a}}^2 + \var(\tns{a})$.
Two RVs $\tns{a}, \tns{b}\in\C{A}$ are \emph{uncorrelated} iff their
covariance vanishes: $\cov(\tns{a},\tns{b}) = 0$, i.e.\ if
their fluctuating parts are orthogonal $\bkt{\tilde{\tns{a}}}{\tilde{\tns{a}}}_{2}=0$.
The two RVs $\tns{a}, \tns{b}\in\C{A}$ are \emph{independent} iff
$\cov(Q_1(\tns{a},\tns{a}^*),Q_2(\tns{b},\tns{b}^*)) = 0$ for all $Q_1, Q_2 \in
\Pi_2^n$ with $\EXP{Q_1(\tns{a},\tns{a}^*)}= \EXP{Q_2(\tns{b},\tns{b}^*)}=0$,
i.e.\ if the centred subspaces of the algebras generated by them are orthogonal,
i.e.\ $(\D{C}\{\tns{a},\tns{a}^*\}\cap\ker\D{E}) \perp
(\D{C}\{\tns{b},\tns{b}^*\}\cap\ker\D{E})$.
In the non-commutative case, the concept of \emph{freeness} and
\emph{free independence} becomes more important, cf.\ \citep{VoiculescuDykemaNica1992,
HiaiPetz2000, MingoSpeicher2017, Speicher2017}, but we shall not further
pursue this topic here.

We have seen that the example algebra $\C{A}_s=\Lp_{0s}(\Omega)$ satisfies
all the requirements and is thus a concrete example of a probability algebra,
and generates the Hilbert space $\Lp_2(\Omega)$,
which is one concrete example of the abstract Hilbert space $\E{V}:=\Lp_2(\C{A})$
for a general probability algebra $\C{A}$.

\subsection{More examples}  \label{SS:xmpls}
For the example algebra
$\C{A}_s=\Lp_{0s}(\Omega)$ it is also well known that one may define
the $\Lp_p$-norms for any $1\le p < \infty$ via $\nd{\phi}_p^p := \EXP{(\phi^* \phi)^{p/2}}
=\int_\Omega |\phi(\omega)|^p\,\D{P}(\di\omega)$.  For $p=\infty$ one
sets $\nd{\phi}_\infty := \text{ess}\sup_\Omega |\phi|$.
The completion of $\C{A}_s=\Lp_{0s}(\Omega)$ in any of the norms $\nd{\cdot}_p$
for $1\le p \le \infty$ gives the familiar Banach spaces $\Lp_p(\Omega)$.
This gives two more concrete examples of probability algebras, namely
$\Lp_\infty(\Omega)$ and $\Lp_{\infty -}(\Omega) := \bigcap_{1\le p <\infty} \Lp_p(\Omega)$.
The last example contains \emph{unbounded} RVs, e.g.\ all the Gaussian RVs.
Obviously one has $\C{A}_s =\Lp_{0s}(\Omega) \subset \Lp_\infty(\Omega) \subset
\Lp_{\infty -}(\Omega) \subset \Lp_0(\Omega)$, i.e.\ the classical simple RVs in $\C{A}_s$ are
a probability sub-algebra of the classical bounded RVs $\Lp_\infty(\Omega)$,
which is a probability sub-algebra
of the algebra $\Lp_{\infty -}(\Omega)$ of unbounded RVs which have finite moments of any order,
which in turn is a sub-$^*$-algebra of the $^*$-algebra of all RVs,
which is not a probability algebra as not every element has a finite expected value.

One more classical example which should be mentioned is the case when
$\Omega$ is in addition a compact Hausdorff topological space,
the $\sigma$-algebra $\F{A}$ is the Borel algebra $\F{B}(\Omega)$, and the
probability measure a Radon measure.  Then the RVs given by the continuous complex-valued
functions $\Ck(\Omega;\D{C})$---for brevity only $\Ck(\Omega)$---are a
sub-probability algebra of $\Lp_\infty(\Omega)$, in fact a $C^*$-algebra---
a Banach space in the $\nd{\cdot}_\infty$ norm such that $\nd{\phi\psi}_\infty
\le \nd{\phi}_\infty\nd{\psi}_\infty$ and $\nd{\phi\phi^*}_\infty=
\nd{\phi}_\infty\nd{\phi^*}_\infty = \nd{\phi}_\infty^2$ such that
the product and adjoint are continuous---called
the \emph{uniform} algebra on $\Omega$.

These are all examples of classical commutative resp.\ Abelian algebras
of RVs with the state the usual Lebesgue integral 
(i.e.\ the usual expected value) w.r.t\ the measure $\D{P}$.
The bounded RVs  $\Lp_\infty(\Omega)$ are
a \emph{maximal Abelian $W^*$-algebra} \citep{Segal1978}
---a $W^*$-algebra is in simplest terms defined as a
$C^*$-algebra which as Banach space is the dual of another Banach space.
It may be shown conversely that any complex
maximal Abelian $W^*$-probability algebra $\C{A}$
is isomorphic to an $\Lp_\infty$-algebra on a probability space,
a result that will be used in the sequel---this is the Segal representation.
Thus the algebraic approach to probability can completely recover the
classical approach due to Kolmogorov which starts from measure spaces
and defines RVs as measurable functions.
Similarly it can be shown that unital Abelian $C^*$-algebras are
isomorphic to the uniform algebra on a compact space---the
Gel'fand representation.
Abelian algebras of this kind are therefore often called `function algebras'.

Let us now consider some non-commutative examples.  A simple one is
$\D{M}(\D{C},n)=\D{C}^{n\times n}$, the algebra of complex
$n\times n$ matrices with complex conjugate transposition as involution.
The language of the algebra is completely the same, except
that projections in the abstract setting---which are self-adjoint--are
called \emph{orthogonal projections} here.  This kind of algebra
corresponds to RVs which can take no more than $n$ different values.
Let $\vek{\vrho}\in\D{M}(\D{C},n)$ be a self-adjoint
positive definite matrix with $\tr \vek{\vrho}=1$, called a \emph{density matrix}.
Then $\D{E}_{\vek{\vrho}}(\vek{A}):= \tr (\vek{\vrho}\vek{A})$ is a faithful state.
Of course any sub-algebra of $\D{M}(\D{C},n)$ which contains the identity
matrix is another example, and  the diagonal matrices
are an example of a commutative sub-algebra.  
More powerful is the algebra $\D{M}(\Lp_{\infty}(\Omega),n)$
of $n\times n$ random matrices with entries from $\Lp_{\infty}(\Omega)$,
and the expectation is the expected value of a matrix state, i.e.\ for
$\tnb{A}\in\D{M}(\Lp_{\infty}(\Omega),n)$ one may set
$\EXP{\tnb{A}} := \int_\Omega \D{E}_{\vek{\vrho}}(\tnb{A}(\omega))\, \D{P}(\di\omega)$.

An example generalising the previous case
is $\E{L}(\C{H})$, the algebra of bounded linear
maps on a complex Hilbert space $\C{H}$ with the adjoint taking the
r\^ole of the involution, or any unital sub-algebra thereof.
$\E{L}(\C{H})$ is a $W^*$-algebra, non-commutative if $\dim\C{H}>1$.
If $\vrho\in\E{L}(\C{H})$ is a nuclear resp.\ trace-class positive
definite operator with unit trace $\tr \vrho = 1$---called again a density matrix---then
a state may be defined for $A\in\E{L}(\C{H})$ as $\D{E}_\vrho(A):=\tr(\vrho A)$.
The example is in some way universal, as
with the \emph{Gel'fand-Naimark-Segal} (GNS) construction any algebra
with faithful state may be embedded (faithfully represented) into an algebra of operators
on a complex Hilbert space \citep{Segal1978,
segal58-TAMS, segalNonlin1969, VoiculescuDykemaNica1992}; namely $\tns{a}\in\C{A}$
is represented as $L_{\tns{a}}: \C{A}\ni \tns{b} \mapsto \tns{ab}
\in \C{A}$ in $\E{L}(\Lp_2(\C{A}))$.

When the Hilbert space $\C{H}$ in question is a Lebesgue space $\Lp_2(\Omega)$,
then any $\kappa\in\Lp_{\infty}(\Omega)$ can be represented as a linear map $M_\kappa:
\Lp_2(\Omega)\ni\vphi \mapsto M_\kappa\vphi=\kappa\vphi\in\Lp_2(\Omega)$.
Thus the Abelian algebra $\Lp_{\infty}(\Omega)$ is represented as a maximal Abelian
$W^*$-sub-algebra of $\E{L}(\Lp_2(\Omega))$, it is called the
\emph{multiplication algebra} of $\Lp_2(\Omega)$.

\subsection{Weights, spectrum, and spectral calculus}  \label{SS:spec-calc}
In this abstract setting we have now seen RVs and their expectation and what
can be deduced from these concepts.  The question arises now as to what
an actual observation or sample of such an RV really is.
To this end a bit more theory is needed.  First it turns out that with
non-commuting observables, in an experiment or other observation,
only commuting observables (self-adjoint elements) can be observed
simultaneously \citep{Whittle2000}.  This is implied by 
the \emph{uncertainty relation}.
Let $\tns{a}, \tns{b} \in\C{A}$ be two self-adjoint elements resp.\ observables,
and $[\tns{a},\tns{b}] = \tns{ab}-\tns{ba}$ be their commutator.
The Cauchy-Bunyakovsky-Schwarz inequality for non-commutative
variables easily gives the \emph{uncertainty relation} 
$\var(a)\var(b) \ge \EXP{\ii [\tns{a},\tns{b}]}^2/4\ge 0$;
where the the expected value on the right hand side is real, as it is
easy to see that $\ii [\tns{a},\tns{b}]$ is self-adjoint.
Once say $\tns{a}$ has been observed, it is known and its variance vanishes.
This shows that it is not possible to observe $\tns{a}$ and $\tns{b}$ simultaneously,
unless they commute.

Therefore the way to approach this is to consider for some observation or
experiment all relevant commuting RVs which can be observed simultaneously,
say $\tns{a}_1, \dots, \tns{a}_k \in \C{A}$.
They, and hence any powers or polynomials in commuting variables of them
can be observed simultaneously, in fact any element of the Abelian sub-probability
algebra $\C{A}_x:=\D{C}[\tns{a}_1,\dots,\tns{a}_k]\subseteq\C{A}$ generated by them.
We shall shortly add more functions beyond polynomials to this list.

As $\tns{a}_1, \dots, \tns{a}_k$ commute, so do the linear operators
$L_{\tns{a}_1},\dots,L_{\tns{a}_k}$ in the GNS-represen\-tation, and
the algebra $\C{L}_x:=\D{C}[L_{\tns{a}_1},\dots,L_{\tns{a}_k}]
\subseteq\E{L}(\Lp_2(\C{A}))$ generated by them is an Abelian algebra
isomorphic to $\C{A}_x$.
It is worthwhile at this point to remember that for linear operators
the fact that they commute means that they have the same spectral resolution,
and the Gel'fand representation of
Abelian $C^*$-algebras and the Segal representation of maximal Abelian $W^*$-algebras
can now be used \citep{gelfand64-vol4, Segal1978, DautrayLions3}.  This can
in fact be employed to obtain a version of the spectral theorem for linear operators.
We defer this for a moment in order to point out the importance
of spectral theory to the subject.

The concept of a state as a self-adjoint positive normalised linear functional
was already introduced.  The set of all possible states $S(\C{A}_x)$
is clearly a subset of the dual $\C{A}_x^*$, and due to the normalisation
they are actually on the unit ball of $\C{A}_x^*$.  One can easily show that
$S(\C{A}_x)$ is a closed, convex, and hence weak-* compact subset of
the unit ball of the dual.  The extreme points of $S(\C{A}_x)$ are
called \emph{pure states}, and their convex combinations are weak-* dense in $S(\C{A}_x)$.
In the case of classical RVs, the states are naturally represented by
probability measures, which are known to form a convex weak-* compact subset
of the unit ball in the space of all measures of bounded total variation.
The extreme points in that case are well known to be Dirac-$\updelta$-measures.

A \emph{weight}, or more specifically a representational weight, also
called a multiplicative character,
$\tns{\alpha}\in S(\C{A}_x)$ is a special kind of state, namely one that is also an
algebra *-homomorphism $\C{A}_x\to\D{C}$.  This means that for $\tns{b},
\tns{c}\in \C{A}_x$ and $\eta, \zeta\in\D{C}$ it holds not only that
$\tns{\alpha}(\eta b + \zeta c)= \ip{\tns{\alpha}}{\eta \tns{b} + 
\zeta \tns{c}} = \eta \tns{\alpha}(\tns{b}) + \zeta \tns{\alpha}(\tns{c})$
(linearity), but also that $\tns{\alpha}(\tns{b}^*) = (\tns{\alpha}(\tns{a}))^*$ and
$\tns{\alpha}(\tns{bc}) = \tns{\alpha}(\tns{b})\tns{\alpha}(\tns{c})$.  The set of all weights
---one-dimensional representations of $\C{A}_x$---is denoted by $\hat{\C{A}_x}$
and is called the \emph{spectrum} of $\C{A}_x$;
it is also a weak-* compact subset $\hat{\C{A}_x}\subset S(\C{A}_x)
\subset B_1(0)\subset \C{A}_x^*$ of the unit ball of the dual.
In the case of classical algebras of RVs the Dirac-$\updelta$-measures are a
good example of weights.

The best known meaning of the term spectrum is certainly when used with regard
to a linear map or an element $\tns{c}\in\C{A}$ as the set $\sigma(\tns{c})=
\{ \lambda\in\D{C} \mid \tns{c} - \lambda \tns{e} \text{ is not invertible} \}$.
Now let $\tns{\alpha}\in\hat{\C{A}_x}$ be any weight, and $\tns{b}\in\C{A}_x$.
If $\tns{b}$ is invertible with inverse $\tns{b}^{-1}$, then $\tns{e}=
\tns{b}\tns{b}^{-1}$ implies $1=\tns{\alpha}(\tns{e})=
\tns{\alpha}(\tns{b}\tns{b}^{-1})=\tns{\alpha}(\tns{b})\tns{\alpha}(\tns{b}^{-1})$,
and hence $\tns{\alpha}(\tns{b})\ne 0$.  Invertible elements can thus not be
mapped to 0 by any weight, i.e.\ any element in the spectrum $\hat{\C{A}_x}$.
Looking at $\tns{b}=\tns{c}-\tns{\alpha}(\tns{c}) \tns{e}$, one sees that
$\tns{\alpha}(\tns{b})= \tns{\alpha}(\tns{c}-\tns{\alpha}(\tns{c}) \tns{e})=
\tns{\alpha}(\tns{c})-\tns{\alpha}(\tns{c})\tns{\alpha}(\tns{e})=0$,
hence $\tns{b}=\tns{c}-\tns{\alpha}(\tns{c}) \tns{e}$
can not be invertible and therefore $\tns{\alpha}(\tns{c})\in\sigma(\tns{c})$ for any
weight $\tns{\alpha}\in\hat{\C{A}_x}$.  This explains the name spectrum for
the set of weights $\hat{\C{A}_x}$, i.e.\ each $\tns{\alpha}(\tns{c})$ is
in the spectrum of $\tns{c}$.  In fact, for any $\lambda \in\sigma(\tns{c})$ there is a
$\tns{\alpha}\in\hat{\C{A}_x}$ such that $\tns{\alpha}(\tns{c})=\lambda$.

The interpretation now is that when one observes a
RV, i.e.\ sees a sample, then one sees the action of some weight on the RV.
Hence the possible values (sample observations) of an abstract RV $\tns{a}\in\C{A}$ are
given by the action of all weights on the RV,
$\{ \tns{\alpha}(\tns{a})= \ip{\tns{\alpha}}{\tns{a}}\mid 
\tns{\alpha}\in\hat{\C{A}_x}\}$.  Therefore one
concludes that all possible observations of a RV $\tns{a}$ are given by
its spectrum $\sigma(\tns{a})$; and as the observables are self-adjoint
the spectrum is real, $\sigma(\tns{a})\subseteq\D{R}$.

Considering general non-commutative probability algebras, the spectrum
of the algebra is often empty as there are no non-zero one-dimensional
representations---another sign that these observables
cannot be observed simultaneously---but in the case of Abelian algebras
like $\C{A}_x$ or $\C{L}_x$, the ones we are considering when examining a concrete experiment
or observation, the Gel'fand and Segal representations tell us that the spectrum
is rich enough.  One may hence use spectral theory of linear
operators to determine the set of possible values, as $\tns{a}\in\C{A}_x$ and
$L_{\tns{a}}\in\C{L}_x$ in the GNS-construction have the same spectrum.

The representation theorems state \citep{Segal1978} that an Abelian probability
algebra is isomorphic to a sub-algebra of $\Lp_\infty(\C{X})$ on a compact
Hausdorff space $\C{X}$.  In fact, the compact space may be chosen as
$\C{X}:=\hat{\C{A}}_x$.  The version of the spectral theorem for linear operators
which is most useful here---and will be used again for a different purpose in
\refS{correlat}---is that an Abelian algebra of operators like $\C{L}_x$ is
not only isomorphic but unitarily equivalent to a sub-algebra of the
multiplication algebra on some measure space $\C{Y}$ \citep{Segal1978, DautrayLions3}
with total measure equal to unity, i.e.\ a classical probability space.
The spectrum of such a multiplication operator $M_\kappa$ with the function or RV
$\kappa\in\Lp_\infty(\C{Y})$ \citep{Segal1978} is the \emph{essential range}
of the function $\kappa$.
Hence any of the commuting RVs $\tns{a}_\ell$ resp.\ $L_{\tns{a}_\ell}$ is
represented by a multiplication operator $M_{\kappa_\ell}$, and hence
as algebra by an RV $\kappa_\ell\in\Lp_\infty(\C{Y})$.  We may thus say
that $\sigma(\tns{a}_\ell)=\sigma(L_{\tns{a}_\ell})=
\sigma(M_{\kappa_\ell})=\sigma(\kappa_\ell)= \text{ess range } \kappa_\ell$.

In the classical framework where RVs are measurable maps on a probability
space, one important and relevant fact is that the composition of measurable
functions is again a measurable function, and one can form new RVs by
applying a measurable function to an existing RV.  In the algebraic framework
presented so far only polynomials---which are kind of natural when dealing
with algebras---have appeared.  Now if $f:\D{R}\to\D{R}$---or more
generally $f:\sigma(\tns{a}_\ell)\subseteq\D{R}\to\D{R}$---is an essentially
bounded measurable function,
so is $\gamma=f\circ\kappa_\ell\in\Lp_\infty(\C{Y})$.  Hence there is a corresponding
$M_\gamma:=f(M_{\kappa_\ell})$ in the multiplication algebra, and a
$L_{\tns{g}}:=f(L_{a_\ell}) \in\E{L}(\Lp_2(\C{A}))$, and a $\tns{g}:=f(a_\ell)$
in the weak-* closure of $\C{A}_x$.  This defines the function $f$ now on the
algebra $\C{L}_x$ or $\C{A}_x$, and is the essence of spectral calculus,
used here to obtain new RVs by applying a measurable function $f$.

\subsection{Extensions}  \label{SS:extensions}
With the spectral calculus in place, one may define non-commutative analogues
of the classical $\Lp_p$-spaces for all $1\le p \le\infty$ by
extending any probability algebra $\C{A}$ through completion in
a certain uniform topology, and not just for $p=2$ as above.  
First note that for a positive element $\tns{a}=\tns{b}^*\tns{b}\in\C{A}$
one can always find a unique positive $\tns{c}\in\C{A}$ such that $\tns{a}
=\tns{c} \tns{c} = \tns{c}^2$ via spectral calculus, as this
$\tns{c}=\tns{a}^{1/2}\in\C{A}$ is the square root.
This allows one to define for any $\tns{a}\in\C{A}$ the absolute value as the
positive element $\ns{\tns{a}}:= (\tns{a}^* \tns{a})^{1/2}\in\C{A}$.
Similarly one may compute the $p$-th power for real $p>0$.  For $1\le p<\infty$
the expression $\nd{\tns{a}}_p^p := \EXP{\ns{\tns{a}}^p}$ defines the $p$-th power of
a norm.  Completion of $\C{A}$ w.r.t.\ any of those norms gives non-commutative
Banach spaces $\Lp_p(\C{A})$, and this agrees for $p=2$ with the previous definition.
It also immediately gives a new algebra $\Lp_{\infty -}(\C{A}):=
\bigcap_{1\le p<\infty}\Lp_p(\C{A})$.

Recalling the spectral calculus from the end of the previous \refSS{spec-calc},
one may now state that $\Lp_p(\C{A})$ contains elements $f(\tns{a})$ for
$\tns{a}\in\C{A}$ and certain measurable functions $f\in\Lp_0(\sigma(\tns{a}))$.
These measurable functions have to be such that in the representation of the Abelian
probability sub-algebra $\D{C}[\tns{a}]$, where $\tns{a}$ is represented by the
multiplication operator $M_\kappa$ on $\Lp_2(\C{Y})$ with $\kappa\in\Lp_\infty(\C{Y})$,
and where $\sigma(\tns{a})=\sigma(\kappa)= \text{ess range } \kappa$, the composite
function satisfies $f\circ\kappa\in\Lp_p(\C{Y})$.

For $p=\infty$ one has to look at the representation of $\tns{a}\in\C{A}$ through
the linear map $L_{\tns{a}}$ in the GNS-construction above and define the
$\nd{\tns{a}}_\infty:=\nd{L_{\tns{a}}}_{op}$ as the operator norm of $L_{\tns{a}}$,
effectively $\nd{\tns{a}}_\infty := \sup_{\tns{a}\ne 0} \nd{\tns{ba}}_2/\nd{\tns{a}}_2$.
One may also define a topology corresponding to the weak operator topology
through the semi-norms $q_{\tns{b},\tns{c}}(\tns{a}) := 
|\bkt{L_{\tns{a}} \tns{b}}{\tns{c}}_{2}|=|\EXP{\tns{c}^* \tns{a} \tns{b}}|$.
Completion of the sub-algebra $\C{A}_\infty := \{ \tns{a}\mid \nd{\tns{a}}_\infty<\infty\}
\subseteq \C{A}$ with finite $\infty$-norm w.r.t.\ the
uniform locally convex topology generated by the semi-norms $q_{\tns{b},\tns{c}}(\cdot)$
gives the probability $W^*$-algebra $\Lp_\infty(\C{A})$.  This
shows that the $\Lp_p$-spaces of non-commutative RVs can be generated
just as in the classical Abelian case.

As already mentioned, the space $\Lp_2(\C{A})$ is a possible candidate for
the space $\E{V}$ appearing in the probabilistic model $S:\E{S}\to\E{V}$.
Other candidates may be generated by the following very general construction:
if $\C{H}$ is a Hilbert space with inner product
$\bkt{\cdot}{\cdot}_0$, and $A$ a possibly unbounded self-adjoint positive
operator in $\C{H}$ with dense domain $\dom A$, one may via spectral calculus
define $A^s$ for any $s>0$ with dense domain $\dom A^s$.
The positive definite sesqui-linear
form given by $\bkt{f}{g}_s:=\bkt{f}{g}_0+\bkt{A^s f}{g}_0$ for $f,g \in\dom A^s$ 
defines an inner product on $\dom A^s$, the completion of which in the associated
topology defines the densely embedded Hilbert space $\C{H}_s\hookrightarrow\C{H}$.
Obviously one also has dense embeddings $\C{H}_s\hookrightarrow\C{H}_t$ for $s> t>0$.
Identifying $\C{H}$ with its dual and denoting the dual of $\C{H}_s$ by $\C{H}_{-s}$,
one obtains Gel'fand triplets \citep{gelfand64-vol4, DautrayLions3} or `sandwiched'
dense embeddings $\C{H}_s\hookrightarrow\C{H} \hookrightarrow\C{H}_{-s}$ 
of Hilbert spaces.   One may even go a step further
and introduce the projective limit $\C{S}=\varprojlim_{s>0} \C{H}_s$, 
depending on $A$ often a nuclear space, which in our case
usually will be a new probability algebra.  The dual construction of
inductive limit $\C{S}^*=\varinjlim_{s>0} \C{H}_{-s}$ then generates
the dual space of generalised objects, like the distributions in the
sense of Sobolev and Schwartz.

It is worthwhile to recall that the
familiar Sobolev-Hilbert spaces $\Hp^s(\D{R}^n)$ are generated in this way
by taking $\C{H}=\Hp^0(\D{R}^n)=\Lp_2(\D{R}^n)$ and $A=-\Updelta + M_{\ns{x}^2}$,
essentially the negative Laplacian added to a multiplication operator.
Then the Schwartz space of rapidly decaying smooth functions $\E{S}(\D{R}^n)$
is the projective limit and additionally an Abelian algebra, and its dual
$\E{S}^\prime(\D{R}^n)$, the inductive limit, is the Schwartz space of
tempered distributions.

The same device can be used here by choosing $\C{H}=\Lp_2(\C{A})$---a space which
is naturally given by the expectation state---and an appropriate
operator $A$; then all the spaces $\C{H}_t, t\in\D{R}$, are possible candidates for $\E{V}$,
and the `regularity' of the RVs in $\E{V}:=\C{H}_t$ can be controlled by the
parameter $t\in\D{R}$.  For $t<0$ these are spaces of `generalised' RVs, only
defined via the duality, similar to the Sobolev-Hilbert spaces with negative exponent.

One possible classical choice for the linear operator $A$ for
$\C{H}=\Lp_2(\Omega)=\Lp_2(\Lp_{\infty-}(\Omega))$ is the following:
denote by $H^{:n:}, n\in\D{N}_0$, the
$n$-th homogeneous chaos \citep{holdenEtAl96, Janson1997} in Wiener's polynomial
chaos decomposition $\C{H}=\overline{\bigoplus}_{n=0}^\infty H^{:n:}$, and define $A$ by
$A\, h := n\, h$ for any $h\in H^{:n:}$; a self-adjoint operator with spectrum
$\sigma(A)=\D{N}_0$, called the number operator.  More examples of Hilbert
spaces of RVs which can be generated in this way may be found in
\citep{holdenEtAl96, Janson1997}, they are all practically
defined with the help of the Wiener-Itô polynomial chaos expansion and are all
possible candidates for the space $\E{V}$.

\subsection{Weak or generalised distributions}  \label{SS:weak-distr}
In any case, this construction of a unital algebra with involution and faithful
state leads to an inner product and Hilbert space $\E{V}$, and the state $\D{E}$ may
be extended as continuous functional onto the whole space $\E{V}$.  This may be used in
the mapping $S:\E{S}\to\E{V}$ in \refS{parametric}.  With the possibility
of also using non-commutative algebras, this approach also allows to
deal with objects such as random matrices, or more generally random fields of
tensors of even order  \citep{hgmRO-1-2018, hgmRO-2-2018}, which is much
more cumbersome in the traditional measure space approach.  Our first example
$\C{A}_s=\Lp_{0s}(\Omega)$ also indicates that the algebraic approach is more
general and can completely recover the measure space approach \citep{Segal1978,
VoiculescuDykemaNica1992, HiaiPetz2000, Speicher2017, MingoSpeicher2017}.
The state takes the place of the usual expectation operator, and
it has all its usual properties.

Nevertheless, even in the general abstract setting of a probability algebra,
it is possible to define a distribution probability measure
or `law' on $\D{R}$ for any non-commutative self-adjoint RV, i.e.\ an observable.
Classically, for a real-valued or self-adjoint RV
$\phi\in\Lp_0(\Omega)$ the law of $\phi$ is the \emph{push-forward}
$\phi_*\D{P}$ of the probability measure $\D{P}$, given for an element
$\C{B}$ of the Borel-$\sigma$-algebra $\F{B}(\D{R})$ by
$\phi_*\D{P}(\C{B}):=\D{P}(\phi^{-1}(\C{B}))$.  

In the abstract setting, for any $\tns{a}\in\C{A}$ one may define the \emph{law}
of $\tns{a}$ as a map $\tau_{\tns{a}}:\Pi_1\to\D{C}$ which assigns
to any polynomial $Q\in\Pi_1$
the number $\tau_{\tns{a}}(Q):=\EXP{Q(\tns{a})}$.
With $\tns{a}\in\C{A}$ self-adjoint, we know that the spectrum is real:
$\sigma(\tns{a})\subseteq\D{R}$.  Let $\C{J}\subset\D{R}$ be a compact
interval which contains the spectrum $\sigma(\tns{a})$.
The polynomials $\Pi_1^r$ with real co-efficients are known to be dense
in $\Ck(\C{J})$ due to the Stone-Weierstrass theorem, and $\tau_{\tns{a}}$
can be shown to be a continuous map,
hence may be extended to all of $\Ck(\C{J})$.  From the Riesz-Markov representation
theorem it now follows that there is a Radon probability measure $\D{P}_{\tns{a}}$
such that $\int_{\C{J}} Q(t)\,\D{P}_{\tns{a}}(\di t) = \tau_{\tns{a}}(Q)$
for any $Q\in\Pi_1^r$, called the
\emph{distribution measure} or \emph{law} of the self-adjoint RV $\tns{a}\in\C{A}$.

This more general approach via a mapping like $S:\E{S}\to\E{V}$ and abstract
probability algebras $\C{A}$ related to $\E{V}$
is also needed in many concrete analytic situations.  As a simple example,
consider, as in \refS{intro} and \refS{parametric}, a RV $\vsigma$
with values in an infinite dimensional Hilbert space $\E{S}$.  For this to
be an `honest' RV, the push-forward distribution
$\vsigma_*\D{P} = \D{P}\circ\vsigma^{-1}$  of the probability measure
$\D{P}$ should be a $\sigma$-additive measure on the Borel sets
$\F{B}(\E{S})$ of $\E{S}$.  It is well known that on a Hilbert space this
is only possible (Sazonov's theorem, cf.\ e.g.\ \citep{bogachev2017, Sullivan2015})
if the correlation $C_\vsigma$ already mentioned in \refS{parametric} is a nuclear
or trace-class operator.  In particular, there is no iso-Gaussian measure---i.e.\
where $C_\vsigma=I$ is the identity, invariant under unitaries---on an
infinite-dimensional Hilbert space; one has to resort to so-called
cylindrical pro-measures (which are not $\sigma$-additive) or enlargements of the
Hilbert space.

The formulations such as with the mapping $S$ or $P$ from above or \refS{parametric}
circumvent all the difficulties mentioned in the previous paragraph
with non-nuclear correlation or covariance operators,
and such an assignment is called a \emph{weak distribution} or
\emph{generalised} RV
\citep{segal58-TAMS, segalNonlin1969, LGross1962, bogachev2017}
resp.\ a \emph{generalised process} \citep{gelfand64-vol4}.  For example
the aforementioned iso-Gaussian weak distribution resp.\ generalised
process---this is also called \emph{white noise} on the Hilbert space $\E{S}$---is
very simply defined:
Pick any complete orthonormal system $\{\vsigma_n\}_n$ in $\E{S}$ and an
infinite sequence of independent identically distributed (iid)
standard Gaussian RVs $\{\zeta_n\}$ (zero mean, unit variance)  as CONS,
and let $\E{H}$ be the Hilbert space generated by them.  Define a linear map
$W:\E{S}\ni \vsigma_n \mapsto \zeta_n \in \E{H}$, and it is clear that its
covariance is $C_W = W^* W = I$, as $W$ is by construction unitary.  Hence
$W$ defines a weak white noise distribution on $\E{S}$.  Other extensions
covered by this use of weak distributions are the cases when the
covariance has continuous spectrum, as often happens for translation
invariant covariance kernels \citep{Matthies_encicl} which are
diagonalised by the Fourier transform \citep{bracewell}.

From all this we conclude that one may define a stochastic model as a weak distribution
on $\E{S}$ via a linear map $S:\E{S}\to\E{V}$, where $\E{V}$ was generated
by a probability algebra $\C{A}$ as described above, and similarly for
$P:\C{V}\to\E{V}$.  For a conventional probability
model we assume that the algebra is Abelian, but the non-commutative case
is useful to model e.g.\ random matrices or tensor fields
\citep{hgmRO-1-2018, hgmRO-2-2018}.  For a dynamical system like the one
mentioned in \refS{intro}, the equality in the equation is to be
understood in a probabilistically weak sense as just described: both
sides of the equation are mapped into the space $\E{V}$, and have to
be equal as elements of that space, i.e.\ in a $\E{V}$-weak sense. First
we spell out the meaning of the map $P$:
\begin{multline} \label{eq:weak-eq-S}
P(\dot{v}(t)) = P(A(\vsigma, \mu;v(t))) + P(f(\vsigma, \mu;t)) \quad \Leftrightarrow \\
   \forall w \in \C{V}: \quad \bkt{\dot{v}(t)}{w}_{\C{V}} = \bkt{A(\vsigma,\mu;v(t))}{w}_{\C{V}}
   + \bkt{f(\vsigma,\mu;t)}{w}_{\C{V}},
\end{multline}
as an element of $\E{V}$, which in detail in $\E{V}$ means
\begin{equation} \label{eq:weak-eq}
  \forall \vphi\in\E{V}: \quad
  \bkt{P(\dot{v}(t))}{\vphi}_{\E{V}} = 
  \bkt{P(A(\vsigma,\mu;v(t)))}{\vphi}_{\E{V}} + 
  \bkt{P(f(\vsigma, \mu;t))}{\vphi}_{\E{V}}.
\end{equation}
This allows one to deal with a much wider range of probabilistic
situations, including white noise as already alluded to, as well as
white noise or a Wiener process in time, as the It\^o-integral can
be understood as a weak stochastic distribution \citep{holdenEtAl96}.
The way \refeq{eq:weak-eq-S} and \refeq{eq:weak-eq} are formulated
also immediately suggests numerical approximations by Galerkin's
method --- called the \emph{stochastic Galerkin} method \citep{matthiesKeese05cmame} ---
using finite dimensional subspaces $\C{V}_n\subseteq\C{V}$
and $\E{V}_m\subseteq\E{V}$.

It may be noted that this whole development is analogous on how generalised functions
or distributions are introduced in the Sobolev-Schwartz framework.  There they
are linear maps from a `nice' space---in fact an algebra---such as $\E{S}(\D{R}^n)$ into
the \emph{algebra} $\D{C}$ with the distinguished state given as
the identity.  Here the generalised probabilistic models on a Hilbert space
$\E{S}$ are linear maps into an \emph{algebra} $\C{A}$ of random variables with
distinguished state $\D{E}$, which again maps the \emph{algebra}  $\C{A}$ into the
algebra $\D{C}$.

%  $Log: rv-algebra.tex,v $
%  Revision 1.2.1.3  2018/09/01 22:26:30  hgm
%  after revision --- really final
%
%  Revision 1.2.1.2  2018/09/01 14:04:42  hgm
%  after revision --- final
%
%  Revision 1.2.1.1  2018/09/01 01:18:40  hgm
%  after revision
%
%  Revision 1.2  2018/07/05 22:01:58  hgm
%  final for arXiv
%
%  Revision 1.1  2018/06/21 23:39:03  hgm
%  first check-in old file from Lin-Param
%
%
%
%
%
%

%%% Local Variables: 
%%% mode: latex
%%% TeX-master: "../18_RV-algebra-model"
%%% End: 

% !TEX root = ../18_RV-algebra-model.tex
% !TEX encoding = UTF-8 Unicode
% RCSID:       $Id: corr-kernel.tex,v 1.3.1.1 2018/09/01 01:22:28 hgm Exp $
% Author:      $Author: hgm $
% Contact:     wire@tu-bs.de
% =================================
%% texfile{
%%  AUTHOR    = "$Author: hgm $",
%%  VERSION   = "$Revision: 1.3.1.1 $",
%%  DATE      = "$Date: 2018/09/01 01:22:28 $",
%%  FILENAME  = "$RCSfile: corr-kernel.tex,v $"}
%
% =================================

\section{Correlation factorisations}  \label{S:correlat}
The  correlation operators $C=R^*R$, $C_\vsigma=S^*S$, and $C_v=P^*P$ have already
been mentioned in \refS{parametric}.  We shall show the development
in terms of the map $R$ defining the parametric variable $r(\mu)$,
for the maps $S$ and $P$ which define the stochastic content, everything
has to be just repeated with different symbols, which we leave for
the reader.  In general,
one may specify \citep{kreeSoize86, hgmRO-1-2018, hgmRO-2-2018}
a densely defined map $C$ in $\C{U}$ through the bilinear form 
\begin{equation}   \label{eq:IX}
   \forall u, v \in \C{U}:\quad \bkt{Cu}{v}_{\C{U}} := \bkt{Ru}{Rv}_{\C{Q}} .
\end{equation}
The map $C=R^* R$, may be called the 
\emph{`correlation'} operator and is by construction self-adjoint
and positive, and if $R$ is continuous so is $C$.
In case the inner product $\bkt{\cdot}{\cdot}_{\C{Q}}$ comes from a
measure $\vpi$ on $\C{M}$, so that for two functions
$\phi$ and $\psi$ on $\C{M}$, one has
\[
  \bkt{\phi}{\psi}_{\C{Q}} := \int_{\C{M}} \phi(\mu) \psi(\mu) \; \vpi(\di \mu),
  \quad \text{such that} \quad
  C = R^* R = \int_{\C{M}} r(\mu) \otimes r(\mu) \; \vpi(\di \mu),
\]
the usual formula for the correlation.
The space $\C{Q}$ may then be taken as $\C{Q}:=\Lp_2(\C{M},\vpi)$.
A special case is when $\vpi$ is a probability measure,
$\vpi(\C{M}) = 1$, as for $\C{M}\gets\Omega$ and $\vpi\gets\D{P}$,
this inspired the term `correlation' operator.
  In terms of the developments in \refS{alg-RV} the Hilbert space $\C{Q}$ would be replaced
by any of the candidates for $\E{V}$ and instead of $C=R^* R$ we
would be investigating $C_\vsigma=S^*S$ or $C_v = P^* P$.

The spectral theorem for operators in a Hilbert space was already used
in \refS{alg-RV}, but here we start in a gentler way.
To make everything as simple as possible to explain the main underlying idea,
assume first that $C$ is a non-singular \emph{trace class} or \emph{nuclear} operator.
This means that it is compact, the spectrum $\sigma(C)$ is a point spectrum, has a
CONS $\{v_m\}_m \subset \C{U}$ consisting of eigenvectors, with
each eigenvalue $\lambda_m\ge \lambda_{m+1}\dots\ge 0$ positive
and counted decreasingly according to
their finite multiplicity, and has finite trace $\tr C = \sum_m \lambda_m < \infty$.
Then a version of the spectral decomposition of $C$ is
\begin{equation}   \label{eq:XIII}
C  = \sum_m \lambda_m (v_m \otimes v_m) .
\end{equation}
Use this CONS to define a new CONS $\{s_m\}_m$ in $\C{Q}$:
 $\lambda_m^{1/2} s_m := R v_m$, to obtain the corresponding
\emph{singular value decomposition} (SVD) of $R$ and $R^*$:
%\ignore{
\begin{multline}   \label{eq:XIV}
R = \sum_m \sqrt{\lambda_m} (s_m \otimes v_m); \quad
R^* = \sum_m \sqrt{\lambda_m} (v_m \otimes s_m); \\
r(\mu) =  \sum_m \sqrt{\lambda_m} \, s_m(\mu) v_m = \sum_m (R^* s_m)(\mu),
\end{multline}
%}
The set $\vsigma(R)=\{\sqrt{\lambda_m}\}_m = \sqrt{\sigma(C)}\subset \D{R}_+$
are the \emph{singular values} of $R$ and $R^*$.
The last relation is
the so-called \emph{\KL{} expansion} or \emph{proper orthogonal decomposition} (POD).
The finite trace condition of $C$ translates into the fact
that $r$ is in $\C{U}\otimes\C{Q}$.
If in that relation the sum is \emph{truncated} at $n\in\D{N}$, i.e.\
\begin{equation}   \label{eq:best-n-term}
r(\mu) \approx r_{\text{ROM}}(\mu) =  \sum_{m=1}^n \sqrt{\lambda_m} \, 
        s_m(\mu) v_m = \sum_{m=1}^n (R^* s_m)(\mu),
\end{equation}
we obtain the \emph{best $n$-term approximation} to $r(\mu)$ in the norm
of $\C{U}$.
Observe that $r$ is linear in the $s_m$.
This means that by choosing the `co-ordinate transformation'
$\C{M}\ni \mu\mapsto (s_1(\mu),\dots,s_m(\mu),\dots)\in\D{R}^{\D{N}}$ one obtains
a \emph{linear / affine} representation where the first co-ordinates are
the most important ones.  For the stochastic cases $C_\vsigma=S^* S$ 
and $C_v=P^*P$ we point out
again as in \refS{alg-RV} that the nuclearity of $C_\vsigma$ resp.\ $C_v$
is necessary for the existence of a measurable map $\vsigma:\Omega\to\E{S}$
resp.\ $v(\vsigma(\cdot),\mu;t):\Omega\to\C{V}$.

Equivalently this means that $S$ resp.\ $P$ has to be a \emph{Hilbert-Schmidt} operator,
e.g.~\citep{DautrayLions3}, a condition which severely restricts
stochastic models.  There is a practical need to consider more general
classes of correlation operators, as already evidenced in the seminal paper
by Karhunen \citep{Karhunen1947, Karhunen1947-e}, where integral
transforms for representations
as in \refeq{eq:best-n-term} were investigated.  This more general view
is for example necessary to consider homogeneous or stationary random fields
or stochastic processes, cf.\ e.g.~\citep{Matthies_encicl}.

One formulation of the spectral decomposition extending \refeq{eq:XIII},
already used implicitly in  \refS{alg-RV}, which does
not require $C$ to be nuclear \citep{DautrayLions3, Segal1978},
nor do $C$ or $R$ have to be continuous, which was used already
in \refS{alg-RV} and has to be applied here to the Abelian
algebra $\D{C}[C]$, is as follows.
The densely defined self-adjoint and positive operator $C:\C{U}\to\C{U}$
is unitarily equivalent with a multiplication operator $M_{\gamma}$
on an appropriate measure space $\C{T}$,
\begin{equation}  \label{eq:ev-mult}
  C = V M_{\gamma} V^*,
\end{equation}
where the unitary map is $V:\Lp_2(\C{T})\to\C{U}$,
and $M_{\gamma}$ multiplies a $\psi\in\Lp_2(\C{T})$ with a
real-valued function $\gamma$; $M_{\gamma}:\psi\mapsto\gamma\psi$.
In case $C$ is bounded, so is $\gamma\in\Lp_\infty(\C{T})$.
As $C$ is positive, $\gamma(t)\ge 0$ for $t\in\C{T}$,
and the essential range of $\gamma$ is the spectrum of $C$.
In  \refS{alg-RV} this was already used for the Abelian algebra $\D{C}[\tns{a}]$
resp.\ $\D{C}[L_{\tns{a}}]$, which says then that any member of
that algebra is unitarily equivalent to a multiplication operator.

As already indicated, via spectral calculus one may define the square root
$M_{\gamma}^{1/2} := M_{\sqrt{\gamma}}$,
and a factorisation similar to $C=R^* R$ is obtained via
$C= (V M_{\sqrt{\gamma}})(V M_{\sqrt{\gamma}})^* =: G^* G$.
From this factorisation and the spectral decomposition \refeq{eq:ev-mult}
follows another singular value decomposition (SVD) of $R$ and $R^*$, which is
\begin{equation}  \label{eq:ev-mult-svd}
R = U M_{\sqrt{\mu}} V^*,\quad R^* = V M_{\sqrt{\mu}} U^*,
\end{equation}
where $U:\Lp_2(\C{T})\to\C{Q}$ is a unitary operator.  Having $M_{\gamma}^{1/2}$
allows us to compute the square root of $C$: $C^{1/2}=V M_{\gamma}^{1/2} V^*$,
and from it the self-adjoint positive definite factorisation $C=C^{1/2}C^{1/2}$.

Consider now an arbitrary factorisation
$C = B^* B$, where $B:\C{U}\to\C{H}$ is a map to a Hilbert space $\C{H}$.
Any two such factorisations $B_1:\C{U}\to\C{H}_1$ and $B_2:\C{U}\to\C{H}_2$ 
with $C=B_1^*B_1=B_2^*B_2$ are \citep{hgmRO-1-2018}
\emph{unitarily equivalent} in that there is a unitary map $X_{21}:\C{H}_1\to\C{H}_2$
such that $B_2 = X_{21} B_1$.  Each such factorisation is also unitarily
equivalent to $R$, i.e.\  there is a unitary $X:\C{H}\to\C{Q}$ such that
$R= X B$.
For finite dimensional spaces, a favourite choice for such a decomposition
of $C$ is the Cholesky factorisation $C = L L^*$, where $B=L^*$ is
represented by an upper triangular matrix.

Let us go back to the situation of \refeq{eq:XIII} and how the SVD of
the factors $R$ \refeq{eq:XIV} in the factorisation $C=R^*R$ was generated.
In the same way a SVD of any of the factorisations just considered may
be generated with left-singular vectors $h_m := B C^{-1} R^* s_m = B C^{-1/2} v_m$,
plus the analogue of \refeq{eq:ev-mult-svd}, i.e.\
\[
B = \sum_m \sqrt{\lambda_m} (h_m \otimes v_m); \quad
B^* = \sum_m \sqrt{\lambda_m} (v_m \otimes h_m); \quad
r =  \sum_m \sqrt{\lambda_m} \, h_m v_m = \sum_m B^* h_m,
\]
and with $W=X^*U$:
\[
B = W M_{\sqrt{\mu}} V^*,\quad R^* = V M_{\sqrt{\mu}} W^*.
\]
The left-singular vectors $h_m$ can now be thought of living on any of the
spaces which appeared in the factorisation, i.e.\ generically $\C{H}$,
for which we have just seen the examples $\C{H}=\Lp_2(\C{T})$
and $\C{H}=\C{U}$ (not necessarily very useful)  \citep{hgmRO-1-2018}.

Instead of $C=B^* B$, one may of course consider
\begin{equation}  \label{eq:ev-mult-Q}
C_{\C{H}} = B B^* = W M_{\gamma} W^*
\end{equation}
on $\C{H}$, which has the same spectrum as $C$---with $C$ nuclear, $C_{\C{H}}$ is also
nuclear---and the whole game can be repeated by looking at the spectral
decompositions of  $C_{\C{H}}$.

When one takes the special case $\C{H}=\C{Q}$ with $C_{\C{Q}} = R R^*$, we see
that $C_{\C{Q}} s_m = \lambda_m s_m$, and $s_m = U V^* v_m$, as
well as $C_{\C{Q}} = U V^* C V U^*$.   
This abstract equation can be spelt out in more analytical
detail for the special case when the inner product on $\C{Q}$ is given
by a measure $\vpi$ on $\C{P}$, as it then becomes
\begin{equation}   \label{eq:XVIII}
   \bkt{C_{\C{Q}}\phi}{\psi}_{\C{Q}} = \bkt{R^* \vphi}{R^* \psi}_{\C{U}} = 
   \iint_{\C{M}\times\C{M}}  \vphi(\mu_1) \vkappa(\mu_1, \mu_2)
    \psi(\mu_2)\; \vpi(\di \mu_1) \vpi(\di \mu_2),
\end{equation}
i.e.\ $C_{\C{Q}}$ is a Fredholm integral operator with kernel
$\vkappa$---on $\C{Q}$ the kernel is in general not reproducing---and
its spectral decomposition $C_{\C{Q}} = \sum_m \lambda_m s_m \otimes s_m$
is nothing but the familiar theorem of Mercer \citep{courant_hilbert}.  Factorisations
of $C_{\C{Q}}$ are then factorisations of the kernel $\vkappa(\mu_1, \mu_2)$
and the corresponding representations of $r(\mu)$ are obtained by integral
transforms \citep{hgmRO-1-2018, hgmRO-2-2018}, as already indicated by
Karhunen in \citep{Karhunen1947, Karhunen1947-e}.
The abstract setting outlined in this section can now be applied to the analysis
of a great number of different situations, see \citep{hgmRO-1-2018} for more detail.

As already indicated, the spectral decomposition \refeq{eq:ev-mult} allows one
to go beyond the requirement that $C$ be nuclear, but in the case of a probability
assignment the push-forward is not a measure any more on $\C{U}$, but it can still
be useful in the computation considering weak distributions.  Another formulation
of the spectral decomposition in the same vein as \refeq{eq:XIII} allows also
to cover the general case \citep{DautrayLions3, gelfand64-vol4}.  
The space $\C{U}=\overline{\bigoplus}_j\C{U}_j$ can be decomposed into a orthogonal direct sum of
invariant subspaces $\C{U}_j$ on each of which the operator has a simple spectrum.
So we may assume for this that the operator has a simple spectrum, otherwise consider each
subspace $\C{U}_j$ in turn.  It turns out
that one can find a so-called \emph{rigged Hilbert space} or \emph{Gel'fand triplet}:
$\C{N}\hookrightarrow\C{U}\hookrightarrow\C{N}^*$ with  $\C{N}$ nuclear and
a densely embedded in $\C{U}$.  The eigenvalue equation for a self-adjoint
operator $C$ can be written in weak form:  for $\lambda\in\sigma(C)$
find $v_\lambda\in\C{U}$ s.t.\ for all $w\in\C{U}$ $\bkt{w}{Cv_\lambda}=
\lambda\bkt{w}{v_\lambda}$, but there may be no $v_\lambda\in\C{U}$ if
$\lambda$ is merely in the spectrum and not also an eigenvalue.
Using duality, this is now weakened to: for $\lambda\in\sigma(C)$
find $v_\lambda\in\C{N}^*$ s.t.\ for all $w\in\C{N}$
$\ip{Cw}{v_\lambda}=\lambda\ip{w}{v_\lambda}$, and it turns out that
one can find such $v_\lambda\in\C{N}^*$, in the larger space $\C{N}^*$.
With this the \refeq{eq:XIII}
may be generalised, where, as the spectrum $\sigma(C)$ may be continuous,
 the sum in general has to be replaced by an  integral w.r.t.\
a measure $\rho$ on $\sigma(C)\subseteq\D{R}$.  As $C=R^*R$, the operator
$C_{\C{Q}}=RR^*$  has the same spectrum, and can be decomposed in a Gel'fand
triplet or rigged Hilbert space $\C{P}\hookrightarrow\C{Q}\hookrightarrow\C{P}^*$
with $s_\lambda\in \C{P}^*$:
\begin{equation}  \label{eq:ev-integr}
  C = \int_{\sigma(C)} \lambda\; v_\lambda \otimes v_\lambda \, \rho(\di \lambda);
  \qquad
    C_{\C{Q}} = \int_{\sigma(C)} \lambda\; s_\lambda \otimes s_\lambda \, \rho(\di \lambda).
\end{equation}
The $s_\lambda\in\C{P}^*$ may be seen as generalised functions, and both decompositions
together in \refeq{eq:ev-integr} allow to write a SVD-like decomposition of $R$ and $R^*$,
corresponding to \refeq{eq:XIV}, and have a representation of $r(\mu)$ in a weak sense
as a \KL{} integral over $\C{P}^*$-generalised functions:
\begin{multline}   \label{eq:XIV-int}
R = \int_{\sigma(C)} \sqrt{\lambda}\;(s_\lambda \otimes v_\lambda)\,\rho(\di \lambda);\quad
R^* = \int_{\sigma(C)} \sqrt{\lambda}\;(v_\lambda \otimes s_\lambda)\,\rho(\di \lambda);\quad\\
r(\mu) =  \int_{\sigma(C)} \sqrt{\lambda}\;s_\lambda(\mu) v_\lambda\,\rho(\di \lambda) = 
          \int_{\sigma(C)} (R^* s_\lambda)(\mu) \,\rho(\di \lambda).
\end{multline}
One familiar and frequent place where this occurs (e.g.\ \citep{Matthies_encicl})
is the classical spectral representation of a stationary stochastic process
\[ q(t) = \int_{\D{R}} \sqrt{S(\omega)} \exp(\ii \omega t)\, Z(\di \omega), \]
where $\sqrt{S(\omega)}$ is the square root of the spectral density---corresponding
to $\sqrt{\lambda}$---and $Z(\di \omega)$ is a random measure with orthogonal
increments and unit variance.  This random measure corresponds to
$v_\lambda\,\rho(\di \lambda)$ in \refeq{eq:XIV-int}, the space $\C{Q}$
corresponds to $\Lp_2(\D{R})$, the space of generalised functions $\C{P}^*$
corresponds to the Schwartz space of tempered distributions $\E{S}^\prime(\D{R})$,
and the generalised eigenfunction $s_\lambda(\mu)$ corresponds to $\exp(\ii \omega t)$,
a generalised eigenfunction of a stationary covariance kernel
which is in $\E{S}^\prime(\D{R})$ but not in $\Lp_2(\D{R})$ \citep{bracewell}.

%  $Log: corr-kernel.tex,v $
%  Revision 1.3.1.1  2018/09/01 01:22:28  hgm
%  after revision
%
%  Revision 1.3  2018/07/06 01:08:31  hgm
%  really final for arXiv
%
%  Revision 1.2  2018/07/05 22:02:08  hgm
%  final for arXiv
%
%  Revision 1.1  2018/06/21 23:39:35  hgm
%  first check-in old file from Lin-Param
%
%
%
%
%
%

%%% Local Variables: 
%%% mode: latex
%%% TeX-master: "../18_RV-algebra-model"
%%% End: 

% !TEX root = ../18_RV-algebra-model.tex
% !TEX encoding = UTF-8 Unicode
% RCSID:       $Id: conclusion.tex,v 1.2.1.1 2018/09/01 01:28:10 hgm Exp $
% Author:      $Author: hgm $
% Contact:     wire@tu-bs.de
% =================================
%% texfile{
%%  AUTHOR    = "$Author: hgm $",
%%  VERSION   = "$Revision: 1.2.1.1 $",
%%  DATE      = "$Date: 2018/09/01 01:28:10 $",
%%  FILENAME  = "$RCSfile: conclusion.tex,v $"}
%
% =================================

\section{Conclusion} \label{S:concl}
Parametric mappings have been analysed together with random
variables with values in infinite dimensional spaces and their generalisations
via an associated linear map, enabling the analysis by using well known techniques for
the analysis of linear mappings.  In the case of stochastic elements
this leads to what is called weak distributions, a generalisation of the
usual concept of a random variable.

In this connection algebras of random variables, the so-called algebraic approach
to probability, leads to a concise description of the generation of
appropriate spaces of random variables, and can naturally be used to
specify randomness on infinite dimensional spaces via weak distributions.
This has as a fundamental building block, next to the algebra of random variables,
a distinguished self-adjoint, positive, and normalised
linear functional called the state, which may be
interpreted as an expectation operator.  It is this setting that turns out to be
conceptually much simpler than the measure-theoretic point of view,
especially in the infinite dimensional setting.  In particular this allows
a natural approach to random matrices and tensor fields, where the random variables
do not necessarily have to commute, and the interesting object is the behaviour
of their spectra, a distinctly analytic and algebraic concept which is
much more complicated to treat with the usual measure-theoretic background.

The associated linear map leads to 
the self-adjoint and positive definite so-called `correlation operator',
as well as its different factorisations.  Different representations
generate different factorisations and thus allow a uniform analysis
of their behaviour via an analysis of linear maps.  It is in particular
the different factorisations, and especially the spectral
decomposition, which lead to suggestions for reduced order models
and their analysis.

Not only does each separated representation define an associated linear map,
but conversely under the restrictive conditions of a nuclear or trace-class
correlation operator each factorisation induces a \KL- or
proper orthogonal decomposition (POD)-like
separated representation.  The extension of this idea to arbitrary non-nuclear
correlations operators is indicated through integral transforms,
exemplified through the use of appropriate spectral
decompositions, either via multiplication operators or as spectral integrals
with rigged Hilbert spaces.  These representations must be classed as
generalised maps or generalised random variables, they can only be considered
in a duality framework in a weak sense.  This can be seen as an analogy to
how normal generalised functions or distributions in the
Sobolev-Schwartz sense are treated as a dual space of very smooth functions,
and in fact the theoretical treatment follows along similar lines.

As this is a very short note touching on many diverse subjects to show their
interconnection, it can naturally only be brief and in many cases just
provides hints which have to be followed further with the references indicated.
The analytic techniques used are `classical' and have been
developed along with the growth of quantum theory in the 1940s.  It is their
combination and uniform view from the point of linear functional analysis
which is novel here.

\ignore{           %%%% BEGIN ignore
\subsection{Test of fonts} \label{SS:testf}
\paragraph{Slanted:} --- serifs
\[ a, \alpha, A, \Phi, \phi, \vphi, \qquad 
     \vek{a}, \vek{\alpha}, \vek{A}, \vek{\Phi}, \vek{\phi}, \vek{\vphi} \]

\paragraph{Slanted:} --- sans serif
\[ \tns{a}, \tns{\alpha}, \tns{A}, \tns{\Phi}, \tns{\phi}, \tns{\vphi}, \qquad
  \tnb{a}, \tnb{\alpha}, \tnb{A}, \tnb{\Phi}, \tnb{\phi}, \tnb{\vphi} \]
   
\paragraph{Other:}  --- no small letters
\[ \C{E}, \C{Q}, \C{R}; \D{E}, \D{Q}, \D{R}; \E{E}, \E{Q}, \E{R}\]

\paragraph{Fraktur:}
\[ \F{e}, \F{q}, \F{r}; \F{E}, \F{Q}, \F{R} \]

\paragraph{Upright:} --- serifs, no small greek letters
\[ \mrm{a}, \mrm{A}, \mrm{\Phi},\qquad \mat{a}, \mat{A}, \mat{\Phi} \]

\paragraph{Upright:} --- sans serif, no small greek letters
\[ \ops{a}, \ops{A}, \ops{\Phi},\qquad \opb{a}, \opb{A}, \opb{\Phi} \]
}           %%%% END ignore

%  $Log: conclusion.tex,v $
%  Revision 1.2.1.1  2018/09/01 01:28:10  hgm
%  after revision
%
%  Revision 1.2  2018/07/05 22:02:20  hgm
%  final for arXiv
%
%  Revision 1.1  2018/06/21 23:40:05  hgm
%  first check-in old file from Lin-Param
%
%
%
%
%
%

%%% Local Variables: 
%%% mode: latex
%%% TeX-master: "../18_RV-algebra-model"
%%% End: 

%\appendix

% ============================================================================
% Bibliography
% The BibTeX files come from my external common bibtex repository
% multiple .bib-files are simply concatenated (without whitespace!)
%\appendix
%\cleardoublepage 

%\begin{thebibliography}{99}

\bibliography{\thebib/jabbrevlong,\thebib/stochastic,\thebib/fuq-new,%
\thebib/fa,\thebib/mat_BU-1-S,\thebib/phys_D,\thebib/num,\thebib/highdim}

%\end{thebibliography}

{ %\color{gray9}
   \tiny
       \texttt{\RCSId} 
   }

%\clearpage
%\printindex

%\clearpage
%\include{publication/informatikberichte}

\end{document}